\documentclass[journal,letterpaper,twocolumn,twoside]{IEEEtran}

\usepackage{times}
\usepackage[margin=1in]{geometry}
\usepackage{amsfonts,amsmath,amssymb,mathrsfs,float,latexsym,amstext}
\usepackage{amsthm}
\usepackage{graphics,graphicx}
\RequirePackage[OT1]{fontenc}
\usepackage{cite}
\RequirePackage%[colorlinks,citecolor=blue,urlcolor=blue]
{hyperref}
\newcommand{\ignore}[1]{} %%% {} empty inside

\allowdisplaybreaks %% allows to break long multiline formulas (use with align, not with aligned!)

\theoremstyle{plain}
\newtheorem{theorem}{Theorem}%[section]

\newtheorem{lemma}{Lemma}%[section]

\theoremstyle{definition}

\newtheorem{remark}{Remark}%[section] 
\newtheorem{example}{Example}%[section]
 
\font\msbmx=msbm10                   % \emptyset should be changed to \varnothing
\font\msbmvii=msbm7                  % in the paper and this will give the proper symbol
\font\msbmv=msbm5
\def\varnothing{\mathchoice{\mbox{\msbmx\char'077}}%
{\mbox{\msbmx\char'077}}{\mbox{\msbmvii\char'077}}{\mbox{\msbmv\char'077}}}%
\newcommand{\mb}[1]{\mathbf{#1}} %% use {\mb{symb}} provides mathbold roman
\newcommand{\Cb}{{\mb{C}}}

\newcommand{\Xb}{{\mb{X}}}

\def\A{{\mb{C}}}
\def\C{{\mb{C}}}

\newcommand{\mc}[1]{\mathcal{#1}} %% use \mc{symb}
\newcommand{\Rc}{{\mc{R}}}
\newcommand{\Rca}{{\bar{\mc{R}}}}

\newcommand{\Nc}{{\mc{N}}}
\newcommand{\Fc}{{\mc{F}}}

\newcommand{\Ac}{{\mc{A}}}

\newcommand{\Lc}{{\mc{L}}}

\newcommand{\wtT}{\widetilde{T}}
\newcommand{\wtX}{\widetilde{X}}
\newcommand{\wtS}{\widetilde{S}}

\usepackage{mathrsfs}
 \newcommand{\PFA}{\mathsf{PFA}}% prob of false alarm
 \newcommand{\pfa}{\mathsf{PFA}}% prob of false alarm
\newcommand{\Pb}{{\mathsf{P}}} %probability
\newcommand{\EV}{{\mathsf{E}}} % expectation
\newcommand{\Eb}{{\mathsf{E}}}
\newcommand{\Hyp}{{\mathsf{H}}} %  hypothesis
\newcommand{\mrm}[1]{\mathrm{#1}}
\newcommand{\drm}{{\mrm{d}}}

\newcommand{\mbb}[1]{\mathbb{#1}} %use \mbb{symb}
\def\One{\mathchoice{\rm 1\mskip-4.2mu l}{\rm 1\mskip-4.2mu l}
{\rm 1\mskip-4.6mu l}{\rm 1\mskip-5.2mu l}}
\newcommand\Ind[1]{{\One_{\{#1\}}}} % indicator -- use \Ind{A}
\newcommand{\Zbb}{\mbb{Z}} %% discrete real line

\newcommand{\class}{{\mbb{C}}(\alpha,\pi)}
\newcommand{\classal}{{\mbb{C}}(\alpha)}

\def\bbr{{\mathbb R}}

\newcommand{\xra}{\xrightarrow} %% use \xla{\text{prob}} matches the length of arrow

\newcommand{\abs}[1]{\left\vert#1\right\vert}
\newcommand{\set}[1]{\left\{#1\right\}}

\newcommand{\brc}[1]{\left(#1\right)}
\newcommand{\brcs}[1]{\left[#1\right]}

\renewcommand{\le}{\leqslant} % AMS le ge
\renewcommand{\ge}{\geqslant}

\usepackage{color}

\begin{document}
%%%%%%%%%%%%%%%%

%%%%
\title{Asymptotic Optimality of Mixture Rules for Detecting Changes in General Stochastic Models
\thanks{The work was supported in part by the Russian Federation 5-100 program, the Russian Federation Ministry of Science and Education Arctic program and the grant 18-19-00452
from the Russian Science Foundation at the Moscow Institute of Physics and Technology.}}

%%%%
\author{Alexander~G.~ Tartakovsky, \IEEEmembership{Senior~Member,~IEEE}  
 \thanks{ A.G. Tartakovsky is a Head of the Space informatics Laboratory at the Moscow Institute of Physics and Technology, Russia
 and Vice President of AGT StatConsult, Los Angeles, California, USA; e-mail: agt@phystech.edu}
\thanks{Manuscript received October 31, 2017; revised May 22, 2018; accepted May 26, 2018.}
\thanks{Copyright (c) 2017 IEEE. Personal use of this material is permitted.  However, permission to use this material for any other purposes must be obtained from the 
IEEE by sending a request to pubs-permissions@ieee.org.}
}

%%%
\markboth{IEEE Transactions on Information Theory,~Vol.~~, No.~~, ~~~2018 (Accepted)}%
{Tartakovsky: Detecting Changes in General Stochastic Models}

\maketitle

%%%
\begin{abstract}
The paper addresses a sequential changepoint detection problem for a general stochastic model, assuming that the observed data may be non-i.i.d. (i.e., dependent and 
non-identically distributed) and the prior distribution of the change point is arbitrary. Tartakovsky and Veeravalli (2005), 
Baron and Tartakovsky (2006), and, more recently, Tartakovsky (2017) developed a general asymptotic theory of changepoint detection for non-i.i.d.\ stochastic models, assuming the 
certain stability of the log-likelihood ratio process, in the case of simple hypotheses when both pre-change and 
post-change models are completely specified. However, in most applications, the post-change distribution is not completely known. 
In the present paper, we generalize previous results to the case of  parametric uncertainty, assuming the parameter of the post-change distribution is unknown. 
We introduce two detection rules based on mixtures -- the Mixture Shiryaev rule and the Mixture Shiryaev--Roberts rule -- and study their asymptotic properties in the Bayesian context. 
In particular, we provide sufficient conditions under which 
these rules are first-order asymptotically optimal, minimizing moments of the delay to detection as the probability of false alarm approaches zero.
\end{abstract}

%%%

%%%%%%%%%%%%%%%%%%%%%%%%%%%%%%%%%%%%%%%%%%%%%%%%%%%%%%%%%%%%%%%%%

%\begin{keyword}[class=AMS]
%\kwd[Primary ]{60F99; 62L10; 62L15}
%\kwd[; secondary ]{60J05; 60G40}
%\end{keyword}

%%%%%%%%%%%%%%%%%%%%%%%%%%%%%%%%%%%%%%%%%%%%%%
%
\section{Introduction} \label{sec:intro}
\IEEEPARstart{S}uppose $X_1,X_2, \dots$ are random variables observed sequentially, which may change statistical properties at an unknown point in time $\nu\in\{0,1, 2, \dots\}$, so that
$X_1,\dots,X_\nu$ are generated by one stochastic model and $X_{\nu+1},  X_{\nu+2}, \dots$ by another model. The value of the change point $\nu$ is unknown and the fact of change 
must be detected as soon as possible controlling for a risk associated with false detections. 

More specifically, let $\Xb^{n}=(X_{1},\dots,X_{n})$ denote a sample of size $n$ and
let  $\{f_{\theta,n}(X_{n}|\Xb^{n-1})\}_{n\ge 1}$ be a sequence of conditional densities of $X_n$ 
given $\Xb^{n-1}$. If $\nu=\infty$, i.e., there is no change, then the parameter $\theta$ is equal to $\theta_0$, so that $f_{\theta,n}(X_{n}|\Xb^{n-1})=f_{\theta_0,n}(X_{n}|\Xb^{n-1})$ for all 
$n \ge 1$. If $\nu=k<\infty$, then $\theta=\theta_1\neq \theta_0$, so that $f_{\theta,n}(X_{n}|\Xb^{n-1})=f_{\theta_0,n}(X_{n}|\Xb^{n-1})$ for $n \le k$ and 
$f_{\theta,n}(X_{n}|\Xb^{n-1})=f_{\theta_1,n}(X_{n}|\Xb^{n-1})$ for $n > k$. 

A sequential detection rule is a stopping time $T$ with respect to an observed sequence $\{X_{n}\}_{n\ge 1}$. That is,  $T$ is an
integer-valued random variable, such that the event $\{T = n\}$, which denotes stopping and taking an action after observing the sample $\Xb^n$, 
belongs to the sigma-algebra  $\Fc_{n}=\sigma(\Xb^n)$ generated by observations $X_1,\dots,X_n$.  
A false alarm is raised when the detection is declared before the change occurs, $T\le \nu$. 
The goal of the quickest changepoint detection problem is to develop a detection rule that stops as soon as possible after the real change occurs under a given risk of false alarms.

In early stages, the work focused on the i.i.d.\ case where $f_{\theta,n}(X_{n}|\Xb^{n-1})= f_\theta(X_n)$, i.e., when  the observations are independent and identically distributed (i.i.d.) 
according to a distribution with density $f_{\theta_0}(X_n)$ in the pre-change mode and with density $ f_{\theta_1}(X_n)$ in the post-change mode.  In the early 1960s,  
Shiryaev~\cite{ShiryaevTPA63} developed a Bayesian 
sequential changepoint detection theory when $\theta_1$ is known. This theory implies that the detection procedure based on thresholding the posterior probability of the 
change being 
active before the current time is strictly optimal,  minimizing the expected delay to detection in the class of procedures with a given weighted probability of false alarm 
if the prior distribution of the change point is geometric. 
At the beginning of the 1970s, Lorden~\cite{lorden-ams71} showed that Page's CUSUM procedure~\cite{page-bka54} is first-order asymptotically optimal in a minimax sense, 
minimizing the maximal expected delay to 
detection in the class of procedures with the prescribed average run length  to false alarm  (ARL2FA) as ARL2FA approaches infinity. In the mid-1980s, Moustakides~\cite{MoustakidesAS86} 
established exact minimaxity of the CUSUM procedure for any value of the ARL2FA. 
Pollak~\cite{PollakAS85} suggested modifying the conventional Shiryaev--Roberts statistic (see \cite{Shiryaev61a,ShiryaevTPA63,Roberts66}) 
by randomizing the initial condition to make it an equalizer.
His version of the Shiryaev--Roberts statistic starts from a random point sampled from the quasi-stationary distribution of the Shiryaev--Roberts statistic. 
He proved that, for a large ARL2FA, this randomized procedure is asymptotically third-order minimax within an additive vanishing term. 
The articles \cite{MoustPolTarSS09,PolunTartakovskyAS09}  indicate that the
Shiryaev--Roberts--Pollak procedure is not exactly minimax for all values of the ARL2FA by showing that a generalized Shiryaev--Roberts procedure that starts from a specially 
designed deterministic point performs slightly better. Shiryaev~\cite{Shiryaev61a,ShiryaevTPA63} was the first who established exact optimality of the Shiryaev--Roberts detection 
procedure in the problem of 
detecting changes occurring at a far time horizon after many re-runs among multi-cyclic procedures with the prescribed mean time between false alarms 
for detecting a change in the drift of the Brownian motion.  
Pollak and Tartakovsky~\cite{PollakTartakovsky-SS09} extended Shiryaev's result to the discrete-time i.i.d. (not necessarily Gaussian) case. 
Third-order asymptotic optimality of generalized Shiryaev--Roberts procedures with random and deterministic head-starts was established in~\cite{tartakovsky-tpa11}.
Another trend related to evaluation of performance of CUSUM and EWMA detection procedures was initiated by the SPC (statistical process control) community
(see, e.g., \cite{bissell-jrssc69,box-book09,crowder-jqt97,hawkins-book98, hawkins-jqt03,montgomery-book08,vandobben-book68,wetherill-book91,woodward-book64,woodall-jqt97,woodall-jqt00}).

In many practical applications, the i.i.d.\ assumption is too restrictive. The observations may be
either non-identically distributed or correlated or both, i.e., non-i.i.d. 
 Lai~\cite{LaiIEEE98} generalized Lorden's asymptotic theory~\cite{lorden-ams71} for the general non-i.i.d.\ case establishing asymptotic optimality of
the CUSUM procedure under very general conditions in the point-wise, minimax, and Bayesian settings. He also suggested a window-limited version of the CUSUM procedure, 
which is computationally less demanding than a conventional CUSUM, but  still preserves asymptotic optimality properties. 
Tartakovsky and Veeravalli~\cite{TartakovskyVeerTVP05}, Baron and Tartakovsky~\cite{BaronTartakovskySA06}, and Tartakovsky~\cite{TartakovskyIEEEIT2017} generalized Shiryaev's Bayesian theory 
for the general non-i.i.d.\ case and for a wide class of prior distributions.  In particular, it was proved that the Shiryaev detection rule is asymptotically optimal -- it minimizes not only the 
expected delay to detection but also higher moments of the detection delay as the weighted probability of  a false alarm vanishes.   
Fuh and Tartakovsky~\cite{FuhTartakovskyIEEEIT2018} specified the results in 
\cite{TartakovskyVeerTVP05,TartakovskyIEEEIT2017} for finite-state hidden Markov models (HMM), finding sufficient conditions under which the Shiryaev and  Shiryaev--Roberts rules 
are first-order asymptotically optimal, 
assuming that both pre-change and post-change distributions are completely specified, i.e., the post-change parameter $\theta_1$ is known. 
Fuh~\cite{Fuh03} proved first-order asymptotic minimaxity of the CUSUM procedure as the ARL2FA goes to infinity. Pergamenchtchikov 
and  Tartakovsky~\cite{PergTarSISP2016} established point-wise and minimax asymptotic optimality properties  of the Shiryaev--Roberts rule  for the general non-i.i.d.\ stochastic model 
in the class of rules with  the prescribed local conditional probability of false alarm (in the given time interval) as well as presented sufficient conditions for ergodic Markov processes.

In a variety of applications, however, a pre-change distribution is known but the post-change distribution is rarely known completely. 
A more realistic situation is parametric uncertainty when the parameter  
$\theta$  of the post-change distribution is unknown since a putative value of $\theta$ is rarely representative. When the post-change parameter is unknown, so that the 
post-change hypothesis 
``$\Hyp_k^\vartheta: \nu=k, \theta=\vartheta$'', $\vartheta \in \Theta$ is composite, and it is desirable to detect quickly a change in a broad range of possible values,
the natural modification of the CUSUM, Shiryaev and Shiryaev--Roberts procedures is based either on maximizing over $\vartheta$ or weighting over
a mixing measure $W(\vartheta)$ the corresponding statistics tuned to $\theta=\vartheta$.  The maximization leads to the generalized likelihood ratio (GLR)-based procedures and 
weighting to mixtures.
 Lorden~\cite{lorden-ams71}  was the first established first-order asymptotic minimaxity of the GLR-CUSUM procedure for the i.i.d.\ exponential families as the ARL2FA goes to infinity 
 (see also Dragalin~\cite{dragalin} for refined results). Siegmund and Yakir~\cite{Siegmund&Yakir_JSPI_08} established third-order asymptotic minimaxity of the randomized 
 mixture Shiryaev--Roberts--Pollak procedure
 for the exponential family with respect to the maximal Kullback--Leibler information.  Lai~\cite{LaiIEEE98} established point-wise and 
minimax asymptotic optimality of the window-limited mixture CUSUM and GLR-CUSUM procedures for general non-i.i.d.\ models. Further detailed overview and references can be found in 
the monographs~\cite{Basseville&Nikiforov-book93,TNB_book2014}.

A variety of applications where sequential changepoint detection is important are discussed, e.g.,  in 
\cite{basseville-automatica88,basseville-automatica98,Basseville&Nikiforov-book93,Bakutetal-book63, 
 Chenetal-amre2011, Kent, 
 masson-book01, montgomery-book08,  Polunchenkoetal-SLJAS2013, Siegmund2013, Tartakovsky-book91, TarVeerASM2004, Tartakovsky&Brown-IEEEAES08, szor, Tartakovsky-Cybersecurity14, 
 Tartakovskyetal-SM06,Tartakovskyetal-IEEESP06, TNB_book2014,Tartakovskyetal_IEEESP2013}.

In this paper, we generalize the asymptotic Bayesian theory, developed in \cite{TartakovskyVeerTVP05,TartakovskyIEEEIT2017} for a simple post-change hypothesis,
to the more important and typical case of the composite post-change hypothesis where the post-change parameter is unknown. We assume that the observations can have a very general structure, i.e., 
can be dependent and non-identically distributed.
The key assumption in the general asymptotic theory is a stability property of the log-likelihood ratio process between the  ``change'' and ``no-change'' hypotheses, which can be formulated in terms of a 
Law of Large Numbers and rates of convergence, e.g.,
as the $r$-complete convergence of the properly normalized log-likelihood ratio and its adaptive version in the vicinity of the true parameter value. 

The rest of the paper is organized as follows.  In Section~\ref{sec:Procedures}, we introduce the mixture Shiryaev and the mixture Shiryaev--Roberts rules. 
In Section~\ref{sec:Problem}, we formulate the asymptotic optimization problems in the class of changepoint detection procedures with the constraint imposed on the weighted probability of false alarm, 
which we address in the following sections. In Section~\ref{sec:AOMS}, we establish the first-order asymptotic optimality of the mixture Shiryaev rule and, in Section~\ref{sec:AOMSR}, 
we study the performance of the mixture Shiryaev--Roberts rule as the weighted probability of false alarm goes to zero. In Section~\ref{sec:AOIR}, we prove asymptotic optimality of the mixture 
Shiryaev and mixture Shiryaev--Roberts rules in a purely Bayesian setup when the cost of delay in change detection approaches zero. 
In Section~\ref{sec:Ex}, we use several examples to illustrate general results. Section~\ref{sec:Remarks} concludes.

%%%___________________________________________
\section{The Shiryaev and Shiryaev--Roberts Mixture Rules} \label{sec:Procedures}

Let $\Pb_\infty$ denote the probability measure corresponding to the sequence of observations $\{X_n\}_{n\ge 1}$  when there is never a change ($\nu=\infty$) and,  
for $k=0,1,\dots$ and $\vartheta\in\Theta$, let $\Pb_{k,\vartheta}$ denote the measure  corresponding to the sequence  $\{X_n\}_{n\ge 1}$ when $\nu=k<\infty$ and $\theta=\vartheta$ 
(i.e., $X_{\nu+1}$ is the first post-change observation), where $\theta\in\Theta$ is a parameter (possibly multidimensional). 
Further, let $p_{k,\vartheta}(\Xb^n)=p(\Xb^n|\nu=k, \theta=\vartheta)$ denote a joint density of the sample $\Xb^n=(X_1,\dots,X_n)$, i.e., density
of the restriction $\Pb_{k,\vartheta}^{(n)}$ of the measure $\Pb_{k,\vartheta}$ to the sigma-algebra $\Fc_n=\sigma(\Xb^n)$ with respect to a non-degenerate sigma-finite measure.
Let $\{g_{n}(X_{n}|\Xb^{n-1})\}_{n\ge 1}$ and $\{f_{\theta,n}(X_{n}|\Xb^{n-1})\}_{n\ge 1}$ be two sequences of conditional densities of $X_n$ 
given $\Xb^{n-1}$. With this notation, the general non-i.i.d.\ changepoint model, which we are interested in, can be written as
\begin{equation}\label{noniidmodel}
\begin{split}
p_{\nu,\theta}(\Xb^n) & = p_\infty(\Xb^n) = \prod_{i=1}^n g_i(X_i|\Xb^{i-1}) ~~ \text{for}~~ \nu \ge n ,
\\
p_{\nu, \theta}(\Xb^n) & =  \prod_{i=1}^{\nu}  g_i(X_i|\Xb^{i-1}) \times \prod_{i=\nu+1}^{n}  f_{\theta,i}(X_i|\Xb^{i-1}) 
\\  
& ~~ \text{for}~~ \nu < n.
\end{split}
\end{equation}
Therefore, $g_n(X_{n}|\Xb^{n-1})$ and $f_{\theta,n}(X_{n}|\Xb^{n-1})$ are the pre-change and post-change conditional densities.
Note that the post-change densities may depend on the change point~$\nu$, i.e., $f_{\theta,n}(X_n|\Xb^{n-1})= f_{\theta,n}^{(\nu)}(X_n|\Xb^{n-1})$ for $n > \nu$. We omit the superscript $\nu$ for brevity. 
While often the pre-change density $g_n$ belongs to the same parametric family as the post-change one $f_{\theta,n}$, i.e., $g_{n}=f_{\theta_0, n}$ for some known value $\theta_0$, 
this is not necessarily the case, so that we consider a more general scenario. %An example in Subsection~\ref{ssec:Ex2} illustrates this more general case.

Let $\EV_{k,\vartheta}$ and $\EV_\infty$ denote expectations under $\Pb_{k,\vartheta}$ and $\Pb_\infty$, respectively.

The likelihood ratio (LR) of the hypothesis ``$\Hyp_k^\vartheta: \nu=k, \theta=\vartheta$'' that the change occurs at $\nu=k$ with the post-change parameter $\theta=\vartheta$ against the no-change hypothesis 
``$\Hyp_\infty: \nu =\infty$'' based on the sample $\Xb^n=(X_1,\dots,X_n)$ is given by the product
\[
LR_{k, n}(\vartheta) = \prod_{i=k+1}^{n}  \frac{f_{\vartheta,i}(X_i|\Xb^{i-1})}{g_i(X_i|\Xb^{i-1})}, \quad n > k
\]
and we set $LR_{k, n}(\vartheta)=1$ for $n \le k$.

Assume that the change point $\nu$ is a random variable independent of the observations with prior distribution 
$\pi_k=\Pb(\nu=k)$, $k=0,1,2,\dots$ with $\pi_k >0$ for $k\in\{0,1,2, \dots\}=\Zbb_+$. We will also assume that a change point may take negative values, which means that the change has 
occurred by the time the observations became available. However, 
the detailed structure of the distribution $\Pb(\nu=k)$ for $k=-1,-2,\dots$ is not important. The only value which matters is the total probability $q=\Pb(\nu \le -1)$ of the change being in effect before the observations
become available. 

Let $\Lc_n(\vartheta) = f_{\vartheta,n}(X_n|\Xb^{n-1})/g_n(X_i|\Xb^{n-1})$. In \cite{TartakovskyVeerTVP05,TartakovskyIEEEIT2017} for detecting a change from  
$\{g_n(X_{n}|\Xb^{n-1})\}$ to $\{f_{\vartheta,n}(X_{n}|\Xb^{n-1})\}$ it was proposed to use the Shiryaev statistic
\begin{equation}\label{S_stat_noniid}
\begin{split}
S_n(\vartheta) & = \frac{1}{\Pb(\nu \ge n)} \brc{q \prod_{i=1}^n \Lc_i (\vartheta)+ \sum_{k=0}^{n-1} \pi_k \prod_{i=k+1}^{n} \Lc_i(\vartheta)} , 
\\
& \quad n \ge 1, \quad S_0(\vartheta)=q/(1-q),
\end{split}
\end{equation}
where $\prod_{i=j}^s \Lc_i =1$ for $s <j$. 

When the value of the parameter is unknown there are two conventional approaches to overcome uncertainty -- either to maximize or average over  $\vartheta$. 
The second approach is usually referred to as 
{\em Mixtures}. To be more specific, introduce a mixing measure $W(\theta)$, $\int_{\Theta} \rm{d} W(\theta) =1$, which can be interpreted as a prior distribution if needed. 
Define  the average (mixed) LR 
\begin{equation}\label{ALR} 
\Lambda_{k,n}^W = \int_\Theta LR_{k, n}(\vartheta)  \, \mrm{d} W(\vartheta), \quad k < n
\end{equation}
and the statistic
\begin{equation}\label{WS_stat}
\begin{split}
S_n^W & =   \int_\Theta S_n(\vartheta) \, \mrm{d} W(\vartheta) 
\\
& = \frac{1}{\Pb(\nu \ge n)} \brc{ q  \Lambda_{0,n}^W + \sum_{k=0}^{n-1} \pi_k \Lambda_{k,n}^W}  , 
\\
& \quad n \ge 1, ~~ S_0^W=q/(1-q) ,
\end{split}
\end{equation}
where $S_n(\vartheta)$ is the Shiryaev statistic tuned to the parameter $\theta=\vartheta$ defined in \eqref{S_stat_noniid}. We will call this statistic the {\em Mixture Shiryaev} (MS) statistic.

In the sequel, we study the MS detection rule  that stops and raises an alarm as soon as  the statistic $S_n^W$ reaches a positive level $A$, i.e., the MS rule is nothing but
the stopping time
\begin{equation}\label{MS_def}
T_A=\inf\set{n \ge 1: S_n^W \ge A},
\end{equation}
where $A>0$ is a threshold controlling for the false alarm risk. In definitions of stopping times we always set $\inf\{\varnothing\}=\infty$.

Another popular statistic for detecting a change from  $\{g_n(X_{n}|\Xb^{n-1})\}$ to $\{f_{\vartheta,n}(X_{n}|\Xb^{n-1})\}$, which has certain optimality properties 
\cite{PollakTartakovsky-SS09,tartakovsky-tpa11,PolunTartakovskyAS09,TNB_book2014}, is the generalized Shiryaev--Roberts (SR) statistic
\begin{equation}\label{SR_stat_noniid}
\begin{split}
R_n(\vartheta) & =\omega LR_{0, n}(\vartheta) + \sum_{k=0}^{n-1} LR_{k, n}(\vartheta) 
\\
& = \omega  \prod_{i=1}^n \Lc_i(\vartheta) + \sum_{k=1}^n \prod_{i=k}^n \Lc_i(\vartheta), \quad n \ge 1
\end{split}
\end{equation}
with a non-negative head-start $R_0(\vartheta) = \omega$, $\omega \ge 0$.
The mixture counterpart, which we will refer to as  the {\em Mixture Shiryaev--Roberts} (MSR) statistic, is
\begin{equation}\label{MSR_stat}
\begin{split}
R_n^W &= \int_\Theta R_n(\vartheta) \, \mrm{d} W(\vartheta),
\\
& = \omega \Lambda_{0,n}^W +  \sum_{k=0}^{n-1} \Lambda_{k,n}^W,  \quad n \ge 1, ~~ R_0^W=\omega ,
\end{split}
\end{equation}
and the corresponding MSR detection rule  is given by the stopping time
\begin{equation}\label{MSR_def}
\wtT_A=\inf\set{n \ge 1: R_n^W \ge A},
\end{equation}
where $A>0$ is a threshold controlling for the false alarm risk. 

In Section~\ref{sec:AOMS}, we show that the MS detection rule $T_A$ is first-order asymptotically optimal, minimizing moments of the stopping time distribution for the low risk of false alarms 
under very general conditions. In Section~\ref{sec:AOMSR}, we establish asymptotic properties of the MSR rule, showing that it is also asymptotically optimal when the prior distribution 
becomes asymptotically flat, but not in general.

%%%______________________________________________________
\section{Asymptotic Problems}\label{sec:Problem}

Let $\Pb^\pi_\theta(\Ac\times \mc{K})=\sum_{k\in \mc{K}}\,\pi_{k} \Pb_{k,\theta}\left(\Ac\right)$  denote the ``weighted'' probability measure and $\Eb^\pi_\theta$ the corresponding expectation.

For $r \ge 1$, $\nu=k \in \Zbb_+$, and $\theta\in\Theta$, introduce the risk associated with the conditional $r$-th moment of the detection delay
\begin{equation}\label{SCrADD}
\Rc^r_{k,\theta}(T)=  \EV_{k, \theta}\left[(T-k)^r\,|\, T> k \right].
\end{equation}
In a Bayesian setting, the average risk associated with the moments of delay to detection  is 
\begin{equation} \label{Riskdef}
\begin{split}
\Rca^r_{\pi,\theta}(T)& : = \Eb^\pi_\theta [ (T-\nu)^r | T> \nu]
\\
&= \frac{{\displaystyle\sum_{k=0}^\infty} \pi_k \Rc^r_{k,\theta}(T)\Pb_\infty( T >k)}{1-\PFA( T)} ,
\end{split}
\end{equation}
where 
\begin{equation} \label{PFAdef}
\pfa(T)=\Pb^\pi_\theta( T \le \nu)= \sum_{k=0}^\infty \pi_k \Pb_\infty( T \le k) 
\end{equation}
is the weighted probability of false alarm (PFA) that corresponds to the risk associated with a false alarm. Note that in \eqref{Riskdef} and \eqref{PFAdef} we used the fact that
$\Pb_{k,\theta}(T \le k) = \Pb_\infty(T\le k)$ since the event $\{T \le k\}$ depends on the observations $X_1,\dots X_k$ generated by the pre-change probability measure $\Pb_\infty$ since
by our convention $X_k$ is the last pre-change observation if $\nu=k$.

In Section~\ref{sec:AOMS}, we are interested  in the Bayesian optimization problem
\begin{equation}\label{sec:PrbfBayes}
\inf_{\{T: \PFA(T) \le \alpha\}}\,\Rca_{\pi, \theta}^r(T)  \quad \text{for all} ~ \theta\in \Theta.
\end{equation}
However, in general this problem is not manageable for every value of the PFA $\alpha\in (0, 1)$. So we will focus on the asymptotic problem assuming that the PFA $\alpha$ approaches zero. 
Specifically, we will be interested in proving that the MS rule is first-order asymptotically optimal, i.e.,
\begin{equation}\label{FOAOdef}
\lim_{\alpha\to0} \frac{\inf_{T\in\class}\Rca_{\pi, \theta}^r(T)}{\Rca_{\pi, \theta}^r(T_A)} =1   \quad \text{for all} ~ \theta\in \Theta,
\end{equation} 
where $\class=\{T: \PFA(T) \le \alpha\}$ is the class of detection rules for which the PFA does not exceed a prescribed number $\alpha \in (0,1)$.
In addition, we will prove that the MS rule is uniformly first-order asymptotically optimal in a sense of minimizing the conditional risk \eqref{SCrADD} for all change point values $\nu=k\in \Zbb_+$, i.e.,
\begin{equation}\label{FOAOunifdef}
\begin{split}
& \lim_{\alpha\to0} \frac{\inf_{T\in\class}\Rc_{k, \theta}^r(T)}{\Rc_{k, \theta}^r(T_A)} =1   
\\
& \quad \text{for all} ~ \theta\in \Theta ~\text{and all}~ k\in \Zbb_+.
\end{split}
\end{equation} 

In Section~\ref{sec:AOIR}, we consider a ``purely'' Bayes problem with the average (integrated) risk, which is the sum of the PFA and the cost of delay proportional to the $r$-th moment of the detection delay and 
prove that the MS rule is asymptotically optimal when the cost of delay to detection approaches $0$.

Asymptotic properties of the MSR rule $\wtT_A$ will be also established.

For a fixed $\theta \in \Theta$, introduce the log-likelihood ratio (LLR) process $\{\lambda_{k, n}(\theta)\}_{n \ge k+1}$ between the hypotheses $\Hyp_{k,\theta}$ ($k=0,1, \dots$) and $\Hyp_\infty$:
$$
\lambda_{k, n} (\theta) = \sum_{j=k+1}^{n}\, \log \frac{f_{\theta,j}(X_{j}|\Xb^{j-1})}{g_{j}(X_{j}|\Xb^{j-1})}, \quad n >k 
$$
($\lambda_{k, n} (\theta) =0$ for $n \le k$).

Let $k\in \Zbb_+$ and $r>0$. We say that a sequence of the normalized LLRs $\{n^{-1}\lambda_{k, n}(\theta)\}_{n\ge 1}$ converges {\em $r-$completely} to a number $I_\theta$ under the 
probability measure $\Pb_{k,\theta}$ as $n \to \infty$ if
\begin{equation}\label{rcompletedef}
\begin{split}
& \sum_{n=1}^\infty n^{r-1} \Pb_{k,\theta} \set{\abs{n^{-1}\lambda_{k, n}(\theta) - I_\theta} > \varepsilon} < \infty 
\\
& \quad \text{for all}~~ \varepsilon > 0,
\end{split}
\end{equation}
and we say that  $\{n^{-1}\lambda_{k, n}(\theta)\}_{n\ge 1}$ converges  to $I_\theta$ {\em uniformly $r-$completely} as $n\to\infty$ if
\begin{equation}\label{unifrcompletedef}
\begin{split}
& \hspace{-5mm}\sum_{n=1}^\infty n^{r-1} \sup_{0 \le k <\infty} \Pb_{k,\theta}\set{\abs{n^{-1}\lambda_{k, n}(\theta) - I_\theta} > \varepsilon} < \infty 
\\
& \quad \text{for all}~~ \varepsilon > 0 .
\end{split}
\end{equation}

Assume that there exists a positive and finite number $I_\theta$ such that the normalized LLR $n^{-1}\lambda_{k,n+k}(\theta)$ converges to $I_\theta$ $r-$completely.
Then it follows from \cite{TartakovskyIEEEIT2017} that when the parameter $\theta$ is known the 
Shiryaev detection rule that raises an alarm at the first time such that the Shiryaev statistic $S_n(\theta)$
 exceeds threshold $(1-\alpha)/\alpha$  is asymptotically (as $\alpha\to 0$) optimal in  class $\class$.

Below we extend this result to the case where $\theta$ is unknown. Specifically, we will show that the MS rule \eqref{MS_def} is asymptotically optimal in problems \eqref{FOAOdef} 
and \eqref{FOAOunifdef} under condition \eqref{rcompletedef} and some other conditions for a large class of priors and all parameter values $\theta\in\Theta$.

%%____________________________________________________________________________________
\section{Asymptotic Optimality of the Mixture Shiryaev Rule}\label{sec:AOMS} 

To study asymptotic optimality we need certain constraints imposed on the prior distribution $\{\pi_k\}$ and on the asymptotic behavior of the decision statistics as the sample size increases (i.e., on the general
stochastic model \eqref{noniidmodel}). 

The following two conditions are imposed on the prior distribution:
\vspace{2mm}

\noindent $\mb{CP} \mb{1}$. {\em  For some} $0 \le \mu <\infty$,
\begin{equation}\label{Prior}
\lim_{n\to\infty}\frac{1}{n}\abs{\log \sum_{k=n+1}^\infty \pi_k} = \mu.
\end{equation}
\noindent $\mb{CP} {\mb 2}$. {\em If $\mu=0$, then in addition}
\begin{equation}\label{Prior1}
\sum_{k=0}^\infty \pi_k |\log\pi_k|^r < \infty \quad  \text{for some} ~ r\ge 1.
\end{equation}
The class of prior distributions satisfying conditions $\mb{CP} \mb{1}$ and $\mb{CP} \mb{2}$  will be denoted by $\Cb(\mu)$.

Note that if $\mu >0$, then the prior distribution has an exponential right tail, in which case, condition \eqref{Prior1} holds automatically. 
If $\mu=0$, the distribution has a heavy tail, i.e., belongs to the model with a vanishing hazard rate.
However, we cannot allow this distribution to have a too heavy tail, which will generate very large time intervals between change points.
This is guaranteed by condition  $\mb{CP} {\mb 2}$. Note that condition  $\mb{CP} {\mb 1}$ excludes light-tail distributions with unbounded hazard rates (e.g., Gaussian-type or Weibull-type
with the shape parameter $\kappa>1$) for which the time-intervals with a change point are very short. 
In this case, prior information dominates information obtained from the observations, the change can be easily detected at early stages, and the asymptotic analysis is impractical.  
 Note also that constraint \eqref{Prior1} is often guaranteed by finiteness of the $r$-th moment, $\sum_{k=0}^\infty k^r \pi_k<\infty$. 

 For $\delta>0$ define $\Gamma_{\delta,\theta}=\{\vartheta\in\Theta\,:\,\vert \vartheta-\theta\vert<\delta\}$. 
 Regarding the general model for the observations \eqref{noniidmodel}, we assume that the following two conditions are satisfied:\\[2mm]
\noindent $\C_{1}$. {\em  There exists a positive and finite number $I_\theta$ such that $n^{-1} \lambda_{k, k+n}(\theta)$ converges to $I_\theta$ in $\Pb_{k,\theta}$-probability and
for any $k\in \Zbb_+$ and $\varepsilon >0$}
\begin{equation}\label{sec:Pmax}
\begin{split}
& \lim_{N\to\infty} \Pb_{k,\theta}\set{\frac{1}{N}\max_{1 \le n \le N} \lambda_{k, k+n}(\theta) \ge (1+\varepsilon) I_\theta} =0 
\\
& \quad \text{for all}~ \theta\in \Theta ;
\end{split}
\end{equation}

\noindent $\C_{2}$. {\em  For any $\varepsilon>0$ there exists $\delta=\delta_{\varepsilon}>0$ such that $W(\Gamma_{\delta,\theta})>0$
and for every $\theta\in\Theta$, for any  $k \in \Zbb_+$, any $\varepsilon>0$, and for some $r\ge 1$}
\begin{equation}\label{rcompLeft}
\begin{split}
& \Upsilon_{k,r}(\varepsilon, \theta) := 
\\
& \sum_{n=1}^\infty n^{r-1} \Pb_{k,\theta}\brc{\frac{1}{n} \inf_{\vartheta\in\Gamma_{\delta,\theta}}\lambda_{k, k+n}(\vartheta) < I_\theta  - \varepsilon} <\infty.
\end{split}
\end{equation}

Note that condition $\C_{1}$  holds whenever $\lambda_{k, k+n}(\theta)/n$ converges almost surely to $I_\theta$ under $\Pb_{\theta,k}$,
\begin{equation}\label{sec:MaRe.1}
\frac{1}{n}\lambda_{k,k+n}(\theta) \xra[n\to\infty]{\Pb_{k, \theta}-\text{a.s.}} I_\theta  \quad \text{for all}~ \theta\in \Theta.
\end{equation}

In order to establish asymptotic optimality we first obtain, under condition $\A_1$, an asymptotic lower bound for moments 
of the detection delay  $\Rca^r_{\pi,\theta}(T)=\Eb^\pi_\theta\brcs{\brc{ T-\nu}^r| T >\nu}$ and $\Rc^r_{k, \theta}= \Eb_{k, \theta}\brcs{\brc{ T-k}^r| T >k}$
of any detection rule $T$ from class $\class$, and then we show that under 
condition $\C_2$ this bound is attained for the MS rule $T_{A}$ when $A=A_\alpha$ is properly selected. 

Asymptotic lower bounds for all positive moments of the detection delay are specified in the following lemma. Condition \eqref{sec:Pmax} (and hence, the a.s.\ convergence 
condition \eqref{sec:MaRe.1}) is sufficient for this purpose. 

%%Lemma
\begin{lemma}\label{Lem:LB}
Let, for some $\mu\ge 0$, the prior distribution belong to class $\Cb(\mu)$. 
Assume that for some positive and finite function $I(\theta)=I_\theta$, $\theta\in\Theta$ condition $\A_1$ holds. 
Then, for all $r>0$ and all $\theta\in\Theta$
\begin{equation}\label{LBinclass}
\liminf_{\alpha\to0} \frac{{\displaystyle\inf_{ T\in\class}} \Rca^r_{\pi,\theta}(T)}{|\log \alpha|^r} \ge  \frac{1}{(I_\theta+\mu)^r}
\end{equation}
and for every $k \in \Zbb_+$,  all $r >0$, and all $\theta\in\Theta$
\begin{equation}\label{LBkinclass}
\liminf_{\alpha\to0} \frac{{\displaystyle\inf_{ T\in\class}}  \Rc^r_{k,\theta}(T)}{|\log \alpha|^r} \ge  \frac{1}{(I_\theta+\mu)^r}.
\end{equation}
\end{lemma}

%%%Proof
\begin{IEEEproof}
The lower bound \eqref{LBinclass} follows from Lemma~1 in Tartakovsky~\cite{TartakovskyIEEEIT2017}. The proof of \eqref{LBkinclass} is a modification and 
generalization of the argument in the proof of Theorem~1 in \cite{PergTarSISP2016} provided in the Appendix.
\end{IEEEproof}

The following lemma provides the upper bound for the PFA of the MS rule.

%%Lemma
\begin{lemma}\label{Lem:PFAMS} 
For all  $A>q/(1-q)$ and any prior distribution of $\nu$, the PFA of the MS rule  $T_{A}$ satisfies the inequality 
\begin{equation}\label{PFAMSineq}
\PFA(T_A) \le 1/(1+A),
\end{equation} 
so that for $\alpha < 1-q$
\begin{equation}\label{Aalpha}
A=A_\alpha=(1-\alpha)/\alpha \quad\text{implies} \quad \PFA(T_{A_\alpha}) \le \alpha. 
\end{equation}
\end{lemma}

%%Proof
\IEEEproof
Clearly,
\[
\PFA(T_A) = \Eb^\pi [\Pb(T_A \le \nu | \Fc_{T_A}); T_A < \infty]. 
\]
Using the Bayes rule and the fact that $\prod_{i=j+1}^n \Lc_i(\theta) =1$ for $j \ge n$, we obtain   
\begin{align*}
 \Pb(\nu=k|\Fc_n)  = 
%\\
%& \frac{\int_\Theta \pi_k \prod_{i=1}^k g_i(X_i|\Xb^{i-1}) \prod_{i=k+1}^{n} f_{\theta, i}(X_i|\Xb^{i-1}) \, \mrm{d} W (\theta)}{\sum_{j=-\infty}^\infty \int_\Theta \pi_j \prod_{i=1}^j g_i(X_i|\Xb^{i-1}) 
%\prod_{i=j+1}^{n} f_{\theta, i} (X_i|\Xb^{i-1}) \, \mrm{d} W (\theta)}
%\\
%& = \frac{\pi_k \int_\Theta \prod_{i=k+1}^{n} \Lc_i(\theta) \, \mrm{d} W (\theta) }{q \int_\Theta \prod_{i=1}^n \Lc_i(\theta) \, \mrm{d} W (\theta) + \sum_{j=0}^{n-1} \pi_j \int_\Theta \prod_{i=j+1}^{n}  \Lc_i(\theta) 
%\, \mrm{d} W (\theta) + \Pb(\nu\ge n)}
  \frac{\pi_k \Lambda_{k,n}^W}{q \Lambda_{0,n}^W + \sum_{j=0}^{n-1} \pi_j \Lambda_{j,n}^W + \Pb(\nu\ge n)},
\end{align*}
so that 
\begin{align*}
\Pb(\nu \ge n | \Fc_n) &= \sum_{k=n}^\infty  \Pb(\nu=k|\Fc_n) 
\\
& =  \frac{\Pb(\nu\ge n)}{q \Lambda_{0,n}^W + \sum_{j=0}^{n-1} \pi_j \Lambda_{j,n}^W + \Pb(\nu\ge n)}
\\
& = \frac{1}{S_n^W + 1}.
\end{align*}
Therefore, taking into account that $S_{T_A}^W \ge A$ on $\{T_A<\infty\}$, we have
\[
\PFA(T_A) = \Eb^\pi [1/ (1+ S_{T_A}^W); T_A< \infty] \le 1/(1+A)
\]
and the inequality \eqref{PFAMSineq} follows. Implication \eqref{Aalpha} is obvious.
\endIEEEproof

The following theorem is the main result in the general non-i.i.d.\ case, which shows that the MS detection rule is asymptotically optimal to the first order under 
mild conditions for the observations and prior distributions. Its proof is given in the Appendix.
 
%%Theorem
\begin{theorem}\label{Th:FOAOgen} 
 Let $r\ge 1$ and let the prior distribution of the change point belong to class $\Cb(\mu)$.  Assume that for some $0<I_\theta<\infty$, $\theta\in\Theta$, right-tail and left-tail conditions
$\A_{1}$ and $\A_{2}$ are satisfied. 

\noindent {\rm \bf (i)} Then, for all $0<m \le r$ and all $\theta\in\Theta$ as $A\to \infty$
\begin{equation} \label{MADDkAOgen}
\Rc^m_{k,\theta}(T_A)  \sim \brc{\frac{\log A}{I_\theta+\mu}}^m \quad \text{for all}~ k \in \Zbb_+
\end{equation}
and
\begin{equation} \label{MADDAOgen}
\Rca^m_{\pi,\theta}(T_A)  \sim \brc{\frac{\log A}{I_\theta+\mu}}^m .
\end{equation}

\noindent {\rm \bf (ii)}   If $A=A_\alpha$ is so selected that $\PFA(T_{A_\alpha}) \le \alpha$ and $\log A_\alpha\sim |\log\alpha|$ as $\alpha\to0$, in particular $A=A_\alpha=(1-\alpha)/\alpha$, 
where $0<\alpha<1-q$, then $T_{A_\alpha}$ is first-order asymptotically optimal as $\alpha\to0$ in class $\class$,  minimizing moments of the detection delay up to order $r$, i.e., 
for all $0<m \le r$ and all $\theta\in\Theta$ as $\alpha\to0$ 
\begin{equation}\label{FOAOmomentskgen}
\begin{split}
 \inf_{T \in \class} \Rc^m_{k,\theta}(T)  & \sim \brc{\frac{|\log\alpha|}{I_\theta+\mu}}^m \sim \Rc^m_{k,\theta}(T_{A_\alpha}) 
 \\
 & \quad \text{for all}~ k \in \Zbb_+
 \end{split}
 \end{equation}
and
\begin{equation}\label{FOAOmomentsgen}
 \inf_{T \in \class} \Rca^m_{\pi,\theta}(T)    \sim \brc{\frac{|\log\alpha|}{I_\theta+\mu}}^m \sim \Rca^m_{\pi,\theta}(T_{A_\alpha}) .
\end{equation}
\end{theorem}

Theorem~\ref{Th:FOAOgen} covers a very wide class of non-i.i.d.\ models for the observations as well as a large class of prior distributions. 
However, condition \eqref{Prior} does 
not include the case where $\mu$ is strictly positive, but may go to zero, $\mu\to 0$. Indeed, as discussed in detail in \cite{TartakovskyIEEEIT2017} the distributions with an exponential right tail 
that satisfy condition \eqref{Prior} with
$\mu>0$ do not converge as $\mu \to 0$ to heavy-tailed distributions for which $\mu=0$. As a result, the assertions  of Theorem~\ref{Th:FOAOgen} do not hold with $\mu=0$ if
$\mu$ approaches $0$ with an arbitrary rate. The rate has to be matched somehow with $\alpha$.
For this reason, we now consider the case where the prior distribution  
$\pi^\alpha=\{\pi_k^\alpha\}$ of the change point depends on the PFA constraint $\alpha$ and becomes ``flat'' when $\alpha$ vanishes. 

In the next lemma, which is analogous to Lemma~\ref{Lem:LB}, we provide asymptotic lower bounds for moments of the detection delay in class $\classal=\mbb{C}(\alpha,\pi^\alpha)$ when the prior distribution  
$\pi^\alpha=\{\pi_k^\alpha\}$ depends on $\alpha$ and $\mu=\mu_\alpha \to 0$ as $\alpha\to0$.

%%Lemma
\begin{lemma}\label{Lem:LB2}
 Let the prior distribution $\pi^\alpha=\{\pi_k^\alpha\}$ of the change point satisfy condition \eqref{Prior} with 
$\mu> 0$ such that $\mu=\mu_\alpha\to 0$ as $\alpha\to 0$.  Assume that for some $0<I_\theta<\infty$, $\theta\in \Theta$, condition $\A_1$ holds. 
Then, for all $r>0$ and $\theta\in\Theta$
\begin{equation}\label{LBinclass2}
\liminf_{\alpha\to0} \frac{{\displaystyle\inf_{ T\in\classal}} \Rca^r_{\pi^\alpha,\theta}(T)}{|\log \alpha|^r} \ge  \frac{1}{I_\theta^r} 
\end{equation}
and 
\begin{equation}\label{LBinclassk2}
\liminf_{\alpha\to0} \frac{{\displaystyle\inf_{ T\in\classal}} \Rc^r_{k,\theta}(T)}{|\log \alpha|^r} \ge  \frac{1}{I_\theta^r} \quad \text{for all}~ k \in \Zbb_+.
\end{equation}
\end{lemma}

\begin{IEEEproof}
The lower bound \eqref{LBinclass2} follows from Lemma~3 in \cite{TartakovskyIEEEIT2017}. 
A proof of the lower bound \eqref{LBinclassk2} is given in the Appendix.
\end{IEEEproof}

Using this lemma, we now establish first-order asymptotic optimality of the MS rule when $\mu=\mu_\alpha$ approaches zero as $\alpha\to0$. To simplify the proof, we 
strengthen condition $\A_2$ in the following uniform version:\\[2mm]
\noindent $\C_{3}$. {\em  For any $\varepsilon>0$ there exists $\delta=\delta_{\varepsilon}>0$ such that $W(\Gamma_{\delta,\theta}) >0$
and for every $\theta\in\Theta$, for any  $\varepsilon>0$, and for some $r\ge 1$}
\begin{equation}\label{rcompunifLeft}
\begin{split}
& \Upsilon_r(\varepsilon,\theta) :=
\\
&\sum_{n=1}^\infty \, n^{r-1} \, \sup_{k \in \Zbb_+} \Pb_{k,\theta}\Big(\frac{1}{n} \inf_{\vartheta\in \Gamma_{\delta,\theta}} \lambda_{k,k+n}(\vartheta) 
\\
&< I_\theta  - \varepsilon \Big)<\infty  .
\end{split}
\end{equation}

%%Theorem
\begin{theorem}\label{Th:FOAOgen2}
Let $r\ge 1$.  Assume that the prior distribution $\pi^\alpha=\{\pi_k^\alpha\}$ of the change point $\nu$ satisfies condition \eqref{Prior} with $\mu=\mu_\alpha\to 0$ 
as $\alpha\to 0$ and that $\mu_\alpha$ approaches zero at such rate that 
\begin{equation}\label{Prior3}
\lim_{\alpha\to 0} \frac{{\sum_{k=0}^\infty  \pi_k^\alpha |\log \pi_k^\alpha|^r}}{|\log \alpha|^r} = 0.
\end{equation} 
Assume that for some $0<I_\theta<\infty$, $\theta\in\Theta$, the right-tail condition $\A_1$ and the uniform left-tail condition $\A_3$ are satisfied.  
If $A=A_\alpha$ is so selected that $\PFA(T_{A_\alpha}) \le \alpha$ and $\log A_\alpha\sim |\log\alpha|$ as $\alpha\to0$, 
in particular $A_\alpha=(1-\alpha)/\alpha$, then the MS rule $T_{A_\alpha}$ is asymptotically optimal as $\alpha\to0$ in class $\classal$, 
minimizing moments of the detection delay up to order $r$: for all $0<m \le r$ and all $\theta\in\Theta$ as $ \alpha \to 0$
\begin{equation} \label{MADDAOgen2}
\begin{split}
 \inf_{T \in \classal} \Rca^m_{\pi^\alpha,\theta}(T)  \sim  \brc{\frac{|\log\alpha|}{I_\theta}}^m 
\sim \Rca^m_{\pi^\alpha,\theta}(T_{A_\alpha}) 
\end{split}
\end{equation}
and for all $k \in \Zbb_+$
\begin{equation} \label{MADDAOgenk2}
\begin{split}
 \inf_{T \in \classal} \Rc^m_{k,\theta}(T) \sim  \brc{\frac{|\log\alpha|}{I_\theta}}^m  \sim\Rc^m_{k,\theta}(T_{A_\alpha}) .
\end{split}
\end{equation}
\end{theorem}

The proof of this theorem is given in the Appendix.

%%%___________________________________________________________
\section{Asymptotic Performance of the Mixture Shiryaev--Roberts Rule}\label{sec:AOMSR}

Consider now the MSR detection rule $\wtT_A$ defined in \eqref{MSR_stat} and \eqref{MSR_def}.

The following lemma shows how  to select threshold $A=A_\alpha$ in the MSR rule to embed it in class $\class$. Write
\begin{align*}
\bar\nu = \sum_{k=0}^\infty k \, \pi_k & =  (1-q) \sum_{k=1}^\infty k \, \Pb( \nu = k | \nu \ge 0) 
\\
&= (1-q) \Eb[\nu | \nu \ge 0].
\end{align*}
Since we are not interested in negative values of $\nu$ we will refer to $\bar\nu$ as the mean of the prior distribution. Recall that $\omega$ ($\omega \ge 0$) is a head-start
of the MSR statistic $R_n^W$ (see \eqref{MSR_stat}).

%%Lemma
\begin{lemma}\label{Lem:PFASR} %%Lemma 4.1
For all $A >0$ and any prior distribution of $\nu$ with finite mean $\bar\nu$, the PFA of the MSR rule $\wtT_A$ satisfies the inequality 
\begin{equation}\label{PFAMSR}
\PFA(\wtT_A)  \le \frac{\omega b + \bar{\nu}} {A} ,
\end{equation}
where $b=\sum_{k=1}^\infty \pi_k$, so that if
\[
A= A_\alpha =   (\omega b+\bar\nu)/\alpha
\]
then $\PFA(\wtT_{A_\alpha}) \le \alpha$, i.e., $\wtT_{A_\alpha} \in \class$.  
\end{lemma}

%%Proof
\IEEEproof
Evidently, $\Eb_\infty[R_n(\vartheta)| \Fc_{n-1}] = 1+ R_{n-1}(\vartheta)$ and hence
\begin{align*}
\Eb_\infty[R_n^W| \Fc_{n-1}] & = \int_\Theta \mrm{d} W(\vartheta) + \int_\Theta R_{n-1}(\vartheta)\, \mrm{d} W(\vartheta) 
\\
& = 1+ R_{n-1}^W.
\end{align*}
So $\{R_n^W-\omega - n\}_{n\ge 1}$ is a zero-mean $(\Pb_\infty,\Fc_n)-$martingale and the MSR statistic $R_n^W$ is  a $(\Pb_\infty,\Fc_n)-$submartingale
with mean  $\Eb_\infty [R_n^W]=\omega+n$. Applying Doob's submartingale inequality, we obtain that for $j=1,2,\dots$
 \[
 \Pb_\infty(\wtT_A \le j) =\Pb_\infty\brc{\max_{1\le i \le j} R_i^W \ge A} \le (\omega+j)/A
 \]
and $\Pb_\infty(\wtT_A \le 0) =0$. Thus,
\begin{align*}
\PFA(\wtT_{A})  & =\sum_{j=1}^\infty \pi_j \Pb_\infty(\wtT_A \le j) 
\\
&\le \frac{\omega \sum_{j=1}^\infty \pi_j +  \sum_{j=1}^\infty j \pi_j} {A} , 
\end{align*}
which proves inequality \eqref{PFAMSR}. Therefore, assuming $\bar{\nu}<\infty$, we obtain that setting 
$A=A_\alpha = (\omega b+\bar\nu)/\alpha$ implies  $\wtT_{A_\alpha} \in \class$ and the proof is complete. 
\endIEEEproof

The following theorem, whose proof is postponed to the Appendix, establishes asymptotic operating characteristics of the MSR rule $\wtT_A$.

%%Theorem
\begin{theorem}\label{Th:AOCMSR} 
 Let $\bar\nu<\infty$ and $0\le \omega < \infty$. Let $r\ge 1$. 
Assume that for some function $0<I_\theta<\infty$, $\theta \in \Theta$, conditions $\A_1$ and $\A_2$ are satisfied.

\noindent {\rm \bf (i)}  Then, for all $0<m \le r$ and $\theta\in\Theta$
\begin{equation} \label{MomentskMSR1}
\lim_{A\to\infty}  \frac{\Rc^m_{k,\theta}(\wtT_A) }{(\log A)^m} =\frac{1}{I_\theta^m} \quad  \text{for all}~ k \in \Zbb_+
\end{equation}
and
\begin{equation} \label{MomentsMSR1}
\lim_{A\to\infty}  \frac{\Rca^m_{\pi,\theta}(\wtT_A)}{(\log A)^m} =\frac{1}{I_\theta^m} .
\end{equation}

\noindent {\rm \bf (ii)}  If $A=A_\alpha$ is so selected that  $\wtT_{A_\alpha}\in \class$ and $\log A_\alpha\sim |\log\alpha|$ as $\alpha\to0$, in particular $A_\alpha=(\omega b+\bar\nu)/\alpha$,
then for all $0<m \le r$ and $\theta\in\Theta$
\begin{equation}\label{MomentskMSR2}
\lim_{\alpha\to0} \frac{\Rc^m_{k,\theta}(\wtT_{A_\alpha})}{|\log \alpha|^m} =\frac{1}{I_\theta^m} \quad \text{for all}~k \in \Zbb_+
\end{equation}
and
\begin{equation}\label{MomentsMSR2}
\lim_{\alpha\to0} \frac{\Rca^m_{\pi,\theta}(\wtT_{A_\alpha})}{|\log \alpha|^m} =\frac{1}{I_\theta^m} .
\end{equation}
\end{theorem}

The next theorem addresses the case where the head-start $\omega=\omega_\alpha$ of the MSR statistic and the mean value $\bar\nu=\bar\nu_\alpha$ of the prior distribution approach infinity 
as $\alpha\to0$ with a certain rate. The proof is given in the Appendix.

%%Theorem
\begin{theorem}\label{Th:AoptMSR}
Assume  that $\omega_\alpha\to\infty$ and $\bar\nu_\alpha\to\infty$ with such rate that the following condition holds:
\begin{equation}\label{Prior4}
\lim_{\alpha\to 0} \frac{{\log (\omega_\alpha+\bar\nu_\alpha)}}{|\log \alpha|} = 0.
\end{equation} 
Assume further that for some $0<I_\theta<\infty$ and $r \ge 1$ conditions $\A_1$ and $\A_3$ are satisfied.  
If threshold $A_\alpha$ is so selected that $\PFA(T_{A_\alpha}) \le \alpha$  and $\log A_\alpha \sim |\log \alpha|$ as $\alpha \to 0$, in particular 
$A_\alpha=(\omega_\alpha b_\alpha+\bar\nu_\alpha)/\alpha$, then for all $0<m \le r$ and $\theta\in\Theta$, as $\alpha\to0$
\begin{equation} \label{MomentsMSRmu0}
\begin{split}
\Rca^m_{\pi^\alpha,\theta}(\wtT_{A_\alpha})\sim \brc{\frac{|\log\alpha|}{I_\theta}}^m 
\sim \inf_{T \in \classal} \Rca^m_{\pi^\alpha,\theta},
\end{split}
\end{equation}
and for all $k\in \Zbb_+$
\begin{equation} \label{MomentskMSRmu0}
\begin{split}
\Rc^m_{k,\theta}(\wtT_{A_\alpha})\sim \brc{\frac{|\log\alpha|}{I_\theta}}^m 
\sim \inf_{T \in \classal} \Rc^m_{k,\theta}(T).
\end{split}
\end{equation}
Therefore, the MSR rule $\wtT_{A_\alpha}$ is asymptotically optimal as $\alpha\to0$ in class $\classal$, minimizing moments of the detection delay up to order $r$. 
\end{theorem}

%%___________________________________________________________________________________
\section{Asymptotic Optimality with Respect to the Integrated Risk}\label{sec:AOIR}

Instead of the constrained optimization problem \eqref{sec:PrbfBayes} consider now the unconstrained, ``purely'' Bayes problem with the loss function
\[
L_{r}(T,\nu) = \Ind{T\le \nu} + c \, (T-\nu)^r \Ind{T > \nu},
\]
where $c>0$ is the cost of delay per unit of time  and  $r \ge 1$. The unknown parameter $\theta$ is now assumed random and the weight function $W(\vartheta)$ is interpreted as the prior distribution of  $\theta$.  
The expected loss (integrated risk) associated with the detection rule~$T$ is given by
\[
\rho_{\pi, W}^{c,r}(T)= \Pb^\pi(T \le \nu)+ c \, \int_\Theta \Eb_\vartheta^\pi[(T-\nu)^+]^r \, \drm W(\vartheta) .
\]
Below we show that the MS rule~$T_A$ with a certain threshold $A=A_{c,r}$ that depends on the cost $c$ is asymptotically optimal, minimizing the integrated risk $\rho_{\pi, W}^{c,r}(T)$ 
over all stopping times as the cost vanishes, $c\to0$. 

Define
\[
\Rc_{\pi,W}^{r}(T) = \int_\Theta  \Rca_{\pi,\vartheta}^r(T) \, \drm W(\vartheta) .
\] 
Observe that, if we ignore
the overshoot, then $\PFA(T_A) \approx 1/(1+A)$ and that using approximation~\eqref{MADDAOgen} we may expect that for a large $A$
\[
\Rc_{\pi,W}^r(T_A) \approx \int_\Theta \brc{\frac{\log A}{I_\vartheta + \mu}}^r \, \drm W(\vartheta) = (\log A)^r D_{\mu,r}
\]
where
\[
D_{\mu,r} =  \int_\Theta \brc{\frac{1}{I_\vartheta + \mu}}^r \, \drm W(\vartheta).
\]
So for large~$A$ the integrated risk of the MS rule is approximately equal to
\begin{align*}
\rho_{\pi, W}^{c,r}(T_A) & = \PFA(T_A) + c \, [1-\PFA(T_A)] \Rc_{\pi,W}^r(T_A)
\\
&\approx 1/A + c \, D_{\mu,r} (\log A)^r := G_{c,r}(A).
\end{align*}
The threshold value $A=A_{c,r}$ that minimizes $G_{c,r}(A)$, $A>0$, is a solution of the equation
\begin{equation} \label{eq:threshold}
r D_{\mu,r} A (\log A)^{r-1}=1/c.
\end{equation}
In particular, for $r=1$ we obtain $A_{c,1}=1/(c D_{\mu,1})$. Thus, it is reasonable to conjecture that threshold~$A_{c,r}$ optimizes the performance of the MS rule for a small $c$, and hence, makes this rule
asymptotically optimal as $c\to0$.

 In the next theorem, whose proof is given in the Appendix, we establish that the MS rule~$T_{A_{c,r}}$ with  threshold~$A_{c,r}$ that satisfies~\eqref{eq:threshold} is indeed asymptotically optimal
as $c\to 0$ under conditions $\A_1$ and $\A_2$ when the set $\Theta$ is compact.

%%%%%Theorem
\begin{theorem} \label{Th:FOasopt_pureBayes} 
Let the prior distribution of the change point belong to class $\Cb(\mu)$.  Assume that for some $0<I_\theta<\infty$, $\theta\in\Theta$, right-tail and left-tail conditions
$\A_{1}$ and $\A_{2}$ are satisfied and that $\Theta$ is a compact set.  Let $A=A_{c,r}$ be the solution of the equation~\eqref{eq:threshold}. Then, as $ c\to 0$,
\begin{equation} \label{RiskAOgenexpflat}
\inf_{T\ge 0} \rho_{\pi, W}^{c,r}(T) \sim D_{\mu,r} \, c \, |\log c|^r  \sim \rho_{\pi, W}^{c,r}(T_{A_{c,r}}) .
\end{equation}
\end{theorem}

Finally, the results analogous to Theorems~\ref{Th:FOAOgen2} and \ref{Th:AoptMSR} in the case where the prior distribution has an exponential tail, i.e.,  $\mu>0$, but $\mu=\mu_c\to 0$ as $c\to0$ 
also hold for the integrated risk. Specifically, let $D_r=D_{\mu=0,r}$, i.e.,
\[
D_{r} =  \int_\Theta \brc{\frac{1}{I_\vartheta^r}} \, \drm W(\vartheta).
\]
Note that the values of the mean of the prior distribution of the change point $\bar\nu =\sum_{j=1}^\infty j \pi_j^c= \bar\nu_c$, 
the head-start  of the MSR statistic  $\omega=\omega_c$, and the value of $b=\sum_{j=1}^\infty \pi_j^c=b_c$ are the functions of the cost $c$.
The following theorem spells out details. The proof is given in the Appendix.

%%Theorem
\begin{theorem}\label{Th:FOaopt_pureBayesflat}
Assume that for some $0<I_\theta<\infty$, $\theta\in\Theta$, right-tail and left-tail conditions $\A_{1}$ and $\A_{3}$ are satisfied and that $\Theta$ is compact. 

\noindent {\bf (i)} If the prior distribution $\pi^c=\{\pi_k^c\}$ satisfies condition \eqref{Prior} with $\mu=\mu_c\to 0$  as $c \to 0$ 
 at such rate that 
\begin{equation}\label{Prior3c}
\lim_{c\to 0} \frac{{\sum_{k=0}^\infty  \pi_k^c |\log \pi_k^c|^r}}{|\log c|^r} = 0 
\end{equation} 
and threshold $A=A_{c,r}$ of the MS rule $T_{A}$ is the solution of the equation 
\begin{equation} \label{threshold1}
r D_{r} A (\log A)^{r-1}=1/c,
\end{equation}
then, as $ c\to 0$,
\begin{equation} \label{MomentsMSmu0Bayesflat} 
\begin{split}
\inf_{T\ge 0} \rho_{\pi^c, W}^{c,r}(T) \sim D_{r} \, c \, |\log c|^r  \sim \rho_{\pi^c, W}^{c,r}(T_{A_{c,r}}) .
\end{split}
\end{equation}
Therefore, the MS rule $T_{A_{c,r}}$ is asymptotically optimal as $c\to0$. 

\noindent {\bf (ii)} If the head-start $\omega_c$ and the mean of the prior distribution $\bar\nu_c$ approach infinity at such rate that 
\begin{equation}\label{Prior4c}
\lim_{c \to 0} \frac{{\log (\omega_c+\bar\nu_c)}}{|\log c|} = 0
\end{equation}
and if $A=A_{c,r}$ of the MSR rule $\wtT_{A}$  is the solution of the equation 
\begin{equation} \label{threshold2}
r D_{r} A (\log A)^{r-1}=(\omega_c b_c +\bar{\nu}_c)/c,
\end{equation}
then, as $ c\to 0$,
\begin{equation} \label{MomentsMSRmu0Bayesflat}
\begin{split}
\inf_{T\ge 0} \rho_{\pi^c,W}^{c,r}(T) \sim D_{r} \, c \, |\log c|^r  \sim \rho_{\pi^c, W}^{c,r}(\wtT_{A_{c,r}}).
\end{split}
\end{equation}
Therefore, the MSR rule $\wtT_{A_{c,r}}$ is asymptotically optimal as $c\to0$. 
\end{theorem}

%%__________________________
\section{Examples}\label{sec:Ex}

%%%%%Remark
\begin{remark} \label{Rem: SufCond}
Obviously, the following condition implies conditions $\C_2$ and $\C_3$:\\[2mm]
\noindent $\C_{4}$. {\em  For any $\varepsilon>0$ there exists $\delta=\delta_{\varepsilon}>0$ such that $W(\Gamma_{\delta,\theta}) >0$. Let the $\Theta \to\bbr_{+}$  function $I(\theta)=I_\theta$ 
be continuous and assume that  for every compact set $\Theta_c\subseteq \Theta$, every  $\varepsilon>0$, and for some $r\ge 1$}
\begin{equation}\label{rcompLeft*}
\begin{split}
&\Upsilon^{*}_r(\varepsilon, \Theta_c):= \sup_{\theta\in\Theta_c} \Upsilon_r(\varepsilon,\theta) = 
\\
&\sum_{n=1}^\infty \, n^{r-1} \, \sup_{\theta \in \Theta_c} \sup_{k \in \Zbb_+} 
\Pb_{k,\theta}\Big(\frac{1}{n} \inf_{\vartheta \in \Gamma_{\delta,\theta}} \lambda_{k,k+n}(\vartheta) 
\\
&< I_\theta  - \varepsilon\Big)<\infty .
\end{split}
\end{equation} 
Hence, it is sufficient for asymptotic optimality of the MS rule 
as well as for asymptotic results related to the MSR rule. Note also that if there exists a continuous $\Theta\times \Theta \to\bbr_{+}$ function $I(\vartheta,\theta)$ 
such that for any  $\varepsilon>0$, any compact  $\Theta_c \subseteq \Theta$ and for some $r\ge 1$
\begin{equation}\label{rcompSup}
\begin{split}
& \Upsilon^{**}_r(\varepsilon, \Theta_c):=
\\
&\sum_{n=1}^\infty \, n^{r-1} \, \sup_{k \in \Zbb_+} \,\sup_{\theta\in \Theta_c} \Pb_{k,\theta}\Big( \sup_{\vartheta \in\Theta_c} \Big\vert\frac{1}{n}\lambda_{k, k+n}(\vartheta) 
\\
& - I(\vartheta,\theta)\Big\vert > \varepsilon\Big) <\infty ,
\end{split}
\end{equation}
then condition $\C_4$, and hence, conditions $\C_2$ and $\C_3$ are satisfied with $I_\theta= I(\theta,\theta)$ since
\begin{align*}
&\Pb_{k,\theta}\brc{\frac{1}{n} \inf_{|\vartheta - \theta|<\delta} \lambda_{k,k+n}(\vartheta) < I_\theta  - \varepsilon} 
\\
& \le \Pb_{k,\theta}\brc{ \sup_{\vartheta \in\Theta_c} \left\vert\frac{1}{n}\lambda_{k, k+n}(\vartheta) - 
I(\vartheta,\theta)\right\vert > \varepsilon} .
\end{align*}        
\end{remark}

As we will see, conditions $\A_4$ and \eqref{rcompSup} are useful in checking of applicability of theorems in particular examples.

%%____________________________________________________________
%\subsection{Detection of a Signal with Unknown Amplitude in Additive Gaussian Noise}\label{ssec:Ex1}
%

\begin{example}[\em Detection of Signals with Unknown Amplitudes in a Multichannel System]\label{Ex1}
Assume there is a multichannel system with $N$ channels (or alternatively an $N$-sensor system) and one is able to observe the output vector 
$X_n=(X_n^1,\dots, X_n^N)$, $n=1,2,\dots$, where the observations in the $i$th channel are  of the  form  
$$
X_{n}^i=\theta_i S_n^i \Ind{n > \nu}  +\xi_{n}^i,\quad n \ge 1.
$$
Here $\theta_i S_n^i$ is a deterministic signal  with an unknown amplitude $\theta_i >0$ that may appear at an unknown time $\nu$ in additive noise $\xi_n^i$. For the sake of 
simplicity, suppose that all signals appear at the same unknown time $\nu$.
Assume that noises $\{\xi_n^i\}_{n \in\Zbb_+}$, $i=1,\dots,N$, are mutually independent  $p^i$-th order Gaussian autoregressive processes AR$(p^i)$, i.e., 
\begin{equation}\label{sec:Ex.1}
\xi_{n}^i = \sum_{j=1}^{p^i} \beta_j^i\xi_{n-j}^i + w_{n}^i, \quad n \ge 1, 
\end{equation}
where   $\{w_n^i\}_{n\ge 1}$ are mutually independent  i.i.d.\ normal $\Nc(0,1)$ sequences and the initial values $\xi_{1-p^i}^i$, $\xi_{2-p^i}^i$, $\dots, \xi_0^i$ are arbitrary random or 
deterministic numbers, in particular we may set zero initial conditions $\xi_{1-p^i}^i=\xi_{2-p^i}^i=\cdots=\xi_0^i=0$. 
The coefficients $\beta_1^i,\dots,\beta_{p^i}^i$ are known and all roots of the equations $z^{p^i} -\beta_1^i z^{p^i-1} - \cdots - \beta_{p^i}^i=0$ are in the interior of the unit circle, so that 
the AR($p^i$) processes are stable. Let $\varphi(x)=(2\pi)^{-1/2}\,e^{-x^2/2}$ denote density of the standard normal distribution. Define the $p_n^i$-th order residual 
\[
\wtX_n^i = X_n^i- \sum_{j=1}^{p_n^i} \beta_j^i X_{n-j}^i, \quad n \ge 1,
\]
where $p_n^i =p^i$ if $n > p^i$ and $p_n^i =n$ if $n \le p^i$.  Write $\theta=(\theta_1,\dots,\theta_N)$ and $\Theta=(0,\infty)\times \cdots \times (0,\infty)$ ($N$ times). 
It is easy to see that the conditional pre-change density is 
\[
g(X_n|\Xb^{n-1})= \prod_{i=1}^N \varphi(\wtX_{n}^i)
\] 
and the post-change density is
\[
f_{\theta}(X_{n}|\Xb^{n-1})=  \prod_{i=1}^N \varphi(\wtX_{n}^i-\theta_i \wtS_{n}^i), \quad \theta\in \Theta,
\]
where $\wtS_n^i =  S_n^i- \sum_{j=1}^{p_n^i} \beta_j^i S_{n-j}^i$. Obviously, due to the independence of the data across channels for all $k \in \Zbb_+$ and $n \ge 1$  the LLR has the form
$$
\lambda_{k, k+n}(\vartheta) = \sum_{i=1}^N \brcs{\vartheta_i  \sum_{j=k+1}^{k+n} \wtS_{j}^i \wtX_{j}^i -\frac{\vartheta_i^2 \sum_{j=k+1}^{k+n} (\wtS_j^i)^2}{2}} .
$$
Under measure $\Pb_{k,\theta}$ the random variables $\{\wtX_{n}^i\}_{n\ge k+1}$ 
are independent Gaussian random variables with mean $\Eb_{k,\theta}[\wtX_n^i]=\theta_i \wtS_n^i$ 
and unit variance, and hence, under  $ \Pb_{k,\theta} $ the normalized LLR can be written as  
\begin{equation}\label{LLRAR1}
\begin{split}
\frac{1}{n} \lambda_{k, k+n}(\vartheta, \theta) & = \sum_{i=1}^N \frac{\vartheta_i \theta_i-\vartheta_i^2/2}{n} \sum_{j=k+1}^{k+n} (\wtS_{j}^i)^2  
\\
& \quad + \frac{1}{n}  \vartheta_i \sum_{j=k+1}^{k+n} \wtS_j^i \eta_j^i ,
\end{split}
\end{equation}
where $\{\eta_j^i\}_{j\ge k+1}$, $i=1,\dots,N$, are mutually independent sequences of i.i.d.\ standard normal random variables. 

Assume that 
\begin{equation}\label{Dr}
\lim_{n\to \infty} \frac{1}{n}  \sup_{k \in \Zbb_+} \sum_{j=k+1}^{k+n} |\wtS_j^i|^2 = Q_i ,
\end{equation}
where  $0<Q_i <\infty$. This is typically the case in most signal processing applications, e.g., for harmonic signals $S_n^i = \sin(\omega_i n + \phi _n^i)$. 
Then for all $k\in\Zbb_+$ and $\theta\in \Theta$
\[
\frac{1}{n} \lambda_{k,k+n}(\theta) \xra[n\to\infty]{ \Pb_{k,\theta} -\text{a.s.}} \sum_{i=1}^N \frac{\theta_i^2 Q_i}{2} =I_{\theta},
\]
so that condition $\C_1$ holds. Furthermore, since all moments of the LLR are finite it is straightforward to show that conditions \eqref{rcompLeft*} and \eqref{rcompSup}, and hence, conditions 
$\C_2$ and $\C_3$  hold for all $r \ge 1$. Indeed, using \eqref{LLRAR1}, we obtain that $I(\vartheta,\theta) =\sum_{i=1}^N(\vartheta_i \theta_i - \vartheta_i^2/2) Q_i$ and for any $\delta >0$
\begin{align*}
&\Pb_{k,\theta}\brc{ \sup_{\vartheta \in [\theta -\delta, \theta+\delta]} \left\vert\frac{1}{n}\lambda_{k, k+n}(\vartheta) - I(\vartheta,\theta)\right\vert > \varepsilon} 
\\
&= 
\Pb_{k,\theta}\brc{|Y_{k, n}(\theta)| > \varepsilon \sqrt{n}},
\end{align*}
where
\[
Y_{k, n}(\theta) = \sum_{i=1}^N \frac{\theta_i}{\sqrt{n}} \sum_{j=k+1}^{k+n} \wtS_j^i \eta_j^i, \quad n \ge 1
\]
is the sequence of normal random variables with mean zero and variance $\sigma_n^2=  n^{-1}\sum_{i=1}^N \theta_i^2 \sum_{j=k+1}^{k+n} (\wtS_j^i)^2$, 
which by \eqref{Dr} is asymptotic to $\sum_{i=1}^N\theta_i^2 Q_i$. 
Thus, for a sufficiently large $n$ there exists $\delta_0>0$ such that  $\sigma^2_n \le \delta_0 + \sum_{i=1}^N\theta_i^2 Q_i$ and we obtain that for all large $n$
\begin{align*}
&\Pb_{k,\theta}\brc{\sup_{\vartheta \in [\theta -\delta, \theta+\delta]} \left\vert\frac{1}{n}\lambda_{k, k+n}(\vartheta) - I(\vartheta,\theta)\right\vert > \varepsilon} 
\\
& \le
 \Pb\brc{|\hat\eta| > \frac{ \delta_0 + \sum_{i=1}^N\theta_i^2 Q_i}{\sigma_n^2} \frac{\varepsilon \sqrt{n}}{\delta_0 + \sum_{i=1}^N\theta_i^2 Q_i}} ,
 \\
 & \le  \Pb\brc{|\hat\eta| > \frac{\varepsilon \sqrt{n}}{\delta_0 + \sum_{i=1}^N\theta_i^2 Q_i}},
\end{align*}
where $\hat\eta\sim \Nc(0,1)$ is a standard normal random variable. Hence, for all $r \ge 1$
\begin{align*}
&\sum_{n=1}^\infty \, n^{r-1} \, \sup_{k \in \Zbb_+} \,\sup_{\theta\in \Theta_c} 
\Pb_{k,\theta}\Big(\sup_{\vartheta \in [\theta -\delta, \theta+\delta]} \Big\vert\frac{1}{n}\lambda_{k, k+n}(\vartheta) 
\\
&- I(\vartheta,\theta)\Big\vert > \varepsilon\Big)< \infty ,
\end{align*}
which implies \eqref{rcompSup} for all $r \ge 1$.

Thus, the MS  rule minimizes as $\alpha\to0$ all positive moments of the detection delay. All asymptotic assertions for the MSR rule presented in Section~\ref{sec:AOMSR} also hold with 
$I_\theta= \sum_{i=1}^N\theta_i^2 Q_i/2$. In particular, by Theorem~\ref{Th:AoptMSR},  the MSR rule is also asymptotically optimal for all $r \ge 1$ if the prior distribution of the 
change point is either heavy-tailed or asymptotically flat.

Since by condition $\C_2$ the MS and MSR procedures are asymptotically optimal for almost arbitrary mixing distribution $W(\theta)$, in this example it is most convenient to 
select the conjugate prior, $W(\theta) = \prod_{i=1}^N F(\theta_i/v_i)$, where $F(y)$ is a standard normal distribution and $v_i>0$, in which case the MS and MSR statistics can be computed explicitly.

Note that this example arises in certain interesting practical applications, as discussed in \cite{TNB_book2014}. For example, surveillance systems (radar, acoustic, EO/IR) typically deal 
with detecting moving and maneuvering targets that appear at unknown times, and  it is necessary to detect a signal from
a randomly appearing target in clutter and noise with the smallest possible delay. In radar applications, often the signal represents a sequence of modulated pulses and 
clutter/noise can be modeled as a Markov Gaussian process or more generally as a $p$-th order Markov process (see, e.g, \cite{Bakutetal-book63,Richards_radar2014}). In underwater detection of 
objects with active sonars, reverberation creates very strong clutter that represents a correlated process in time \cite{Marage_sonar2013}, so that again the problem can be reduced 
to detection of a signal with an unknown intensity in correlated clutter. In applications related to detection of point and slightly extended objects with EO/IR sensors (on moving and still platforms such as 
space-based, airborne, ship-board, ground-based), sequences of images usually contain a cluttered background which is correlated in space and time, and it is a challenge to detect and track
 weak objects in correlated clutter~\cite{Tartakovsky&Brown-IEEEAES08}. 

Yet another challenging application area where the multichannel model is useful is cyber-security
\cite{Tartakovsky-Cybersecurity14,Tartakovskyetal-SM06,Tartakovskyetal_IEEESP2013}. Malicious intrusion attempts in computer networks (spam campaigns, personal data theft, 
worms, distributed denial-of-service (DDoS) attacks, etc.)  incur significant financial damage and are a severe harm to the integrity of personal information.  
It is therefore essential to devise automated techniques to detect computer  network intrusions as quickly as possible so that an appropriate response can be provided and the 
negative consequences for the users are eliminated. In particular, DDoS attacks typically involve many traffic streams resulting in a large number of packets aimed at 
congesting the target's server or network.  As a result, these attacks usually lead to abrupt changes in network traffic and can be detected by noticing a change in the average
number of packets sent through the victim's link per unit time. Figure~\ref{fig:UDPattack} illustrates how the multichannel anomaly Intrusion Detection System works for detecting a 
real UDP packet storm.  The multichannel MSR algorithm with the AR$(1)$ model and uniform prior $W(\theta_i)$ on a finite interval $[1,5]$ was used.
The first plot shows packet rate. It is seen that there is a slight change in the mean, which is barely visible. 
The second plot shows the behavior of the multi-cyclic MSR statistic $W_n=\log R_n^W$, which is restarted from scratch every time a threshold exceedance occurs. 
Threshold exceedances before the UDP DDoS attack starts (i.e., false alarms) are shown by green dots and the true detections are marked by red dots. 

 \begin{figure*}[t]
 \centerline {\includegraphics[width=0.8\textwidth]{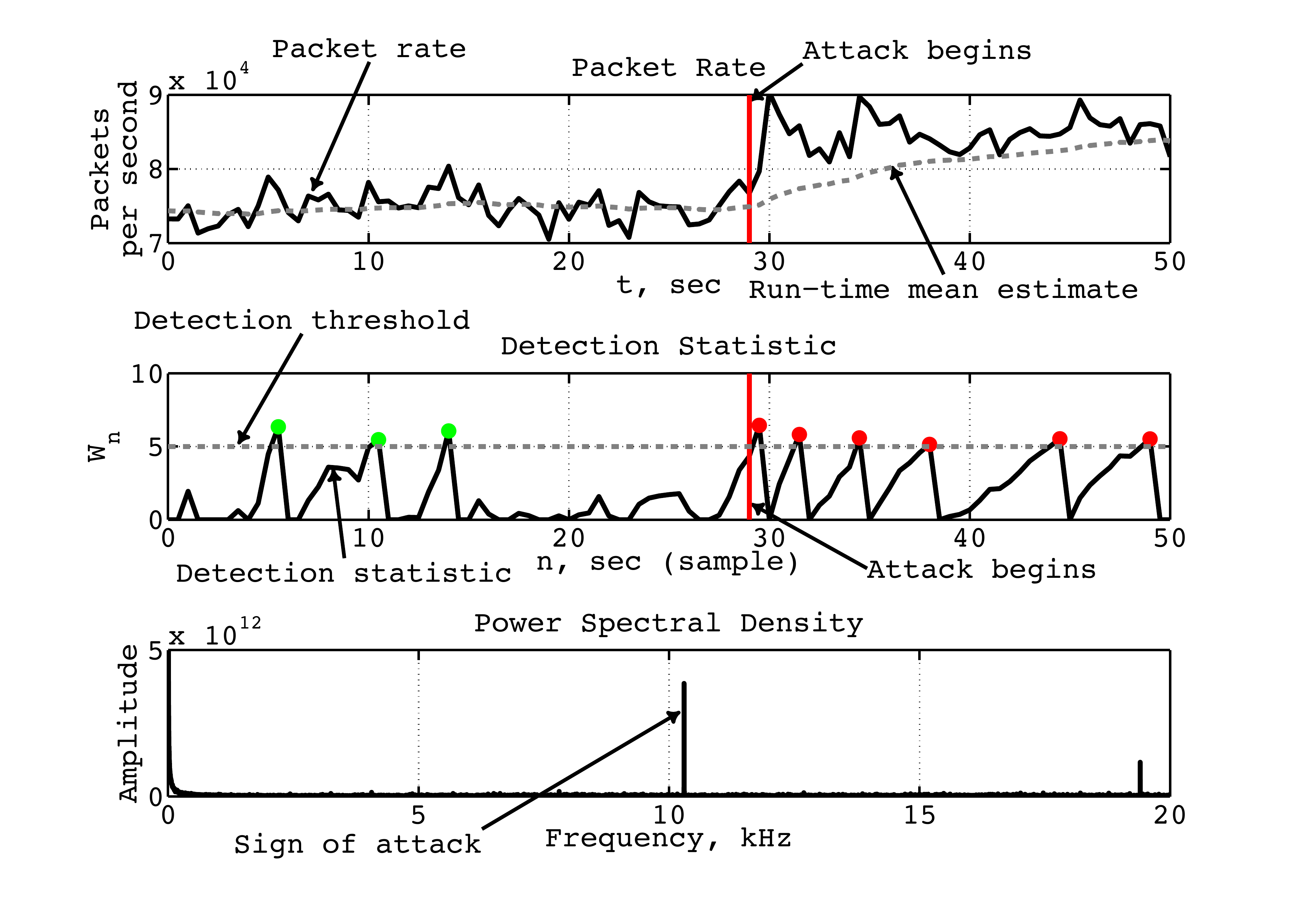}}
 \caption{Detection of the UDP DDoS packet storm attack: upper picture --- raw data (packet rate); bottom  --- log of the MSR statistic.}
 \label{fig:UDPattack}
\end{figure*}
\end{example}

%%______________________________________________________________________________________
%\subsection{Example 3: Detection of Changes in a  Hidden Markov Model} \label{ssec:Ex3}

\begin{example}[\em Detection of Changes in a  Hidden Markov Model] \label{Ex2}

The following example, which deals with a two-state hidden Markov model with i.i.d.\ observations in each state may be of interest, in particular,  for rapid detection and tracking of
sudden spurts and downfalls in activity profiles of terrorist groups that  could be caused by various factors such as changes in the organizational dynamics of terrorist groups, counterterrorism
activity, changing socio-economic and political contexts, etc. In \cite{Raghavanetal-AoAS2013}, based on the analysis of  real data from {\em Fuerzas Armadas Revolucionarias de Colombia}
(FARC) terrorist group from Colombia (RDWTI) it was shown that two-state HMMs can be 
recommended for detecting and tracking 
sudden changes in activity profiles of terrorist groups.
The HMM framework  provides good explanation capability of past/future activity
across a large set of terrorist groups with different ideological attributes.

 Specifically, let $\upsilon_n\in\{1,2\}$ be a two-state Markov chain with the transition matrix   
\begin{align*}
[P_\theta(\upsilon_{n-1}=i, \upsilon_n=l)] = \left[
\begin{array}{cc}
1 - \beta_{\theta} &  \beta_{\theta} \\
\gamma_{\theta} & 1 - \gamma_{\theta}
\end{array}
\right] 
\end{align*}
and stationary initial distribution $\Pb_{\theta}(\upsilon_0 = 2)  = 1-\Pb_{\theta}(\upsilon_0 = 1)= \pi_\theta(2)=\gamma_{\theta}/(\beta_{\theta} + \gamma_{\theta})$ for some 
$\beta_{\theta} , \gamma_{\theta}\in [0 , 1]$, where  the parameter $\theta$ equals $\theta_0$ in the pre-change mode ($\theta_0$ is known) and  $\theta \in \Theta$ in the post-change mode (unknown).  
Suppose that conditioned on $\upsilon_n$  the observations $X_n$  are i.i.d.\ with densities $p_\theta(X_n| \upsilon_n=l) =p_\theta^{(l)}(X_n)$, $l=1,2$.

Introduce the probabilities $P_{\theta, n}  :=  \Pb_{\theta} ( \Xb^n, \upsilon_n = 2)$ and 
$\widetilde{P}_{\theta, n}  :=  \Pb_{\theta} ( \Xb^n, \upsilon_n = 1)$. Straightforward computations show that for $n \ge 1$
\begin{align*}
P_{\theta, n} & =   \left[ P_{\theta, n-1} \, (1-\gamma_\theta) + \widetilde{P}_{\theta, n-1} \, \beta_\theta \right] \, p_{\theta}^{(2)}(X_n) ; %\label{Rec1}
\\
\widetilde{P}_{\theta, n}  & =  \left[ P_{\theta, n-1} \, \gamma_\theta + \widetilde{P}_{\theta, n-1} \, (1-\beta_\theta) \right] \, p_{\theta}^{(1)}(X_n) %\label{Rec2}
\end{align*}
with initial values $P_{\theta, 0}  = \pi_\theta(2)$ and $\widetilde{P}_{\theta,0} =\pi_\theta(1) =1-\pi_\theta(2)$. Denote $p_{0, \theta}(\Xb^n)=p_{\theta} (\Xb^n)$ and $p_{\infty}(\Xb^n)=p_{\theta_0} (\Xb^n)$. Since $p_{\theta} (\Xb^n)=P_{\theta, n}+\widetilde{P}_{\theta,n}$ we obtain that
\[
p_{k,\theta}(\Xb_{k+1}^n | \Xb_0^k) = \frac{P_{\theta, n} + \widetilde{P}_{\theta, n}}{P_{\theta, k} + \widetilde{P}_{\theta, k}}
\]
and
\[
\lambda_{k, k+n}(\theta) = \log \brc{\frac{P_{\theta, k+n} + \widetilde{P}_{\theta, k+n}} {P_{\theta_0, k+n} + \widetilde{P}_{\theta_0, k+n}} 
\cdot \frac{P_{\theta_0, k} + \widetilde{P}_{\theta_0, k}} {P_{\theta, k} + \widetilde{P}_{\theta, k}}}. 
\]

For the sake of simplicity, consider now the symmetric case where $\beta_\theta=\gamma_\theta=1/2$ for all $\theta\in \Theta + \theta_0$. Then
\begin{align*}
P_{\theta,n} & = \frac{1}{2^n} p_\theta^{(2)} (X_n) \prod_{i=1}^{n-1} \brcs{p_\theta^{(1)} (X_n)+ p_\theta^{(2)} (X_n)},
\\
\widetilde{P}_{\theta,n} & = \frac{1}{2^n} p_\theta^{(1)} (X_n) \prod_{i=1}^{n-1} \brcs{p_\theta^{(1)} (X_n)+ p_\theta^{(2)} (X_n)}
\end{align*}
and we obtain that the LLR is
\[
\lambda_{k, k+n}(\theta)  = \sum_{i=k+1}^{k+n} \log \brc{\frac{p_\theta^{(1)} (X_i)+ p_\theta^{(2)} (X_i)}{p_{\theta_0}^{(1)} (X_i)+ p_{\theta_0}^{(2)} (X_i)}}.
\]

Condition $\A_1$ holds with 
\[
I_\theta = \int \log \brc{\frac{p_\theta^{(1)} (x)+ p_\theta^{(2)} (x)}{p_{\theta_0}^{(1)} (x)+ p_{\theta_0}^{(2)} (x)}} \frac{p_\theta^{(1)} (x)+ p_\theta^{(2)} (x)}{2} \, \mrm{d} x 
\]
being the Kullback--Leibler information number since by the SLLN $n^{-1} \lambda_{k, k+n}(\theta) \to I_\theta$ $\Pb_{k,\theta}$-a.s. (assuming that $I_\theta <\infty$). 
Condition \eqref{rcompLeft*} usually holds if  the $(r+1)$-th absolute moment of the increment of the LLR is finite:
\begin{align*}
& \int \abs{\log \brc{\frac{p_\theta^{(1)} (x)+ p_\theta^{(2)} (x)}{p_{\theta_0}^{(1)} (x)+ p_{\theta_0}^{(2)} (x)}}}^{r+1} \times
\\
&
\frac{1}{2} \brcs{p_\theta^{(1)} (x)+ p_\theta^{(2)} (x)} \, \mrm{d} x  < \infty.
\end{align*}

This is the case, for example, if the observations are Gaussian with unit variance and different mean values in 
pre- and post-change modes as well as for different states, i.e.,
$p_\theta^{(l)}(y) = \varphi(y-\mu_\theta^{(l)})$ ($\theta=\theta_0$ or $\theta \neq \theta_0$, $l=1,2$). It is easily verified  that the Kullback--Leibler number  
$I_\theta$ is finite and that condition \eqref{rcompLeft*} is satisfied for all $r \ge 1$. Therefore, in this case, the 
MS and MSR detection rules are asymptotically optimal, minimizing asymptotically all positive moments of the detection delay.
\end{example}

\section{Concluding Remarks}\label{sec:Remarks}

1.  In the case where the increments $\{\Delta \lambda_i(\theta)\}$ of the LLR 
$\lambda_{k,n}(\theta)=\sum_{i=k+1}^n \Delta \lambda_i(\theta)$ are independent (but not necessarily identically distributed), condition $\C_2$ in Theorem~\ref{Th:FOAOgen} and
Theorem~\ref{Th:AOCMSR} and condition $\C_3$ in Theorem~\ref{Th:FOAOgen2} and Theorem~\ref{Th:AoptMSR} can be relaxed in the following condition:
for all $\ell \ge k$ and $k \in \Zbb_+$
\begin{equation}\label{probLeft}
\Pb_{k,\theta}\brc{\frac{1}{n} \int_{\Gamma_{\delta,\theta}} \lambda_{\ell, \ell+n}(\vartheta) \, \mrm{d} W(\vartheta) < I_\theta  - \varepsilon} \xra[n\to\infty]{} 0 .
\end{equation}
More specifically, the MS rule asymptotically minimizes all moments of the delay to detection under the right-tail and left-tail conditions $\A_1$ and \eqref{probLeft}. Also, assertions of 
Theorem~\ref{Th:FOAOgen2} and Theorem~\ref{Th:AoptMSR} for the MSR rule hold for all $r \ge 1$ under conditions $\A_1$ and \eqref{probLeft}. 
The proof can be built by considering the cycles 
$[k+(n-1) N_A, k+n N_A]$, $n =1,2,\dots$ of the length  $N_A = 1+\lfloor \log A/(I_\theta +\mu-\varepsilon) \rfloor$ and slightly modifying of the technique developed by 
Tartakovsky~\cite{TartakovskyIEEEIT2017} in the case of complete knowledge of the post-change distribution.

2. Since we do not assume a class of models for the observations such as Gaussian, Markov or HMM and build the decision statistics on the LLR process 
$\lambda_{k, k+n}(\theta)$,  
it is natural to impose conditions on the behavior of $\lambda_{k, k+n}(\theta)$, which is expressed by conditions $\C_1$, $\C_2$ and $\C_3$, related to the law of large numbers for the LLR and 
rates of convergence in the law of large numbers. The assertions of Theorems~\ref{Th:FOAOgen}--\ref{Th:AoptMSR} hold if 
$n^{-1}\lambda_{k, k+n}(\theta)$ and $n^{-1} \log \Lambda_{k, k+n}^W$ converge uniformly $r$-completely to $I_\theta$ under $\Pb_{k,\theta}$, i.e., when for all $\varepsilon >0$ and $\theta\in\Theta$
\begin{align}
& \sum_{n=1}^\infty \, n^{r-1} \, \sup_{k \in \Zbb_+} \Pb_{k,\theta}\brc{\abs{\frac{1}{n} \lambda_{k, k+n}(\theta)-I_\theta} > \varepsilon} <\infty , \nonumber
\\
&\sum_{n=1}^\infty \, n^{r-1} \, \sup_{k \in \Zbb_+} \Pb_{k,\theta}\brc{\abs{\frac{1}{n}  \log \Lambda_{k,k+n}^W-I_\theta} > \varepsilon} <\infty . \label{rcompA}
\end{align}
However, verifying the $r$-complete convergence condition \eqref{rcompA} for the weighted LLR $\log \Lambda_{k,k+n}^W$ is typically much more difficult than checking 
conditions $\C_2$ and $\C_3$ for the local values of the LLR in the vicinity of the true parameter value. For the simple post-change hypothesis, sufficient conditions for the class of 
ergodic Markov models 
are given in \cite{PergTarSISP2016} and for HMMs in \cite{FuhTartakovskyIEEEIT2018}. In these cases, the LLR process is a Markov random walk and one of the key conditions is finiteness of the 
$(r+1)$-th moment of its increment.

3. As expected, the results indicate that the MSR rule is not asymptotically optimal when the prior distribution of the change point has an exponential tail (i.e., $\mu>0$), 
but it is asymptotically optimal for heavy-tailed prior distributions (i.e., $\mu=0$) and also when $\mu\to0$ with a certain rate.

4. The results show that first-order asymptotic optimality properties of the MS and MSR procedures hold for practically arbitrary weight function $W(\theta)$, 
in particular for any prior that has strictly positive
values on $\Theta$. Therefore, the selection of $W(\theta)$ can be based solely on the computational aspects. The conjugate prior is typically the best choice when possible. However, 
if the parameter is vector and the parameter space is intricate, constructing mixture statistics may be difficult. In this case, discretizing 
the parameter space and selecting the prior $W(\theta=\theta_i)$ concentrated on discrete points $\theta_i$, $i=1,\dots,N$, suggested and discussed in \cite{FellourisTartakovsky-SS2013} 
for the hypothesis testing problems, is perhaps the best option. 
Then one can easily compute the MS and MSR statistics (as long as the LR $\Lambda_{k,k+n}(\theta)$ can be computed) at the expense of losing optimality between the points $\theta_i$ since
 the resulting discrete versions are asymptotically optimal only at the points $\theta_i$.

%%_________________________________________
\section*{Acknowledgement}

The author would like to thank Prof. Sergey Pergamenchtchikov for pointing out that conditions \eqref{rcompLeft*} and \eqref{rcompSup} are sufficient for condition $\C_3$, 
which was useful in verification of $\C_3$ in examples. Thanks also go to referees whose comments improved the presentation.

%%_________________________________________________________
%\appendix

\renewcommand{\theequation}{A.\arabic{equation}}
\setcounter{equation}{0}

%%_________________________________________________________
\section*{Appendix: Proofs}

%%Proof
\begin{IEEEproof}[Proof of Lemma~\ref{Lem:LB}]
For $\varepsilon\in(0,1)$ and $\delta>0$, define $N_{\alpha}=N_{\alpha}(\varepsilon,\delta,\theta)= (1-\varepsilon) |\log\alpha|/(I_\theta+\mu+\delta)$. By the Chebyshev inequality,
\begin{align*}
& \Rc_{k, \theta}^r(T) \ge \Eb_{k, \theta}[( T-k)^+]^r
\\
&\ge N_{\alpha}^r \Pb_{k, \theta}(T  -k > N_{\alpha})
\\
& \ge  N_{\alpha}^r\brcs{\Pb_{k,\theta}( T>k)-\Pb_{k,\theta}(k <  T < k+ N_{\alpha})} ,
\end{align*}
where $\Pb_{k,\theta}( T>k)=\Pb_\infty (T>k)$, so that 
\begin{equation} \label{A1}
\begin{split}
&\inf_{ T\in\class}\Rc_{k,\theta}^r(T) 
\\
&\ge N_{\alpha}^r \Big [ \inf_{T\in\class}\Pb_{\infty}( T>k) 
\\
&- \sup_{ T\in \class}\Pb_{k,\theta}(k <  T < k+N_{\alpha}) \Big].
\end{split}
\end{equation}
Thus, to prove the lower bound \eqref{LBinclass} we need to show that, for arbitrary small $\varepsilon$ and $\delta$ and all fixed $k \in \Zbb_+$, 
\begin{equation}\label{Pinf}
\lim_{\alpha\to0} \inf_{T\in\class}\Pb_{\infty}( T>k) =1
\end{equation}
and
\begin{equation}\label{Psupclasszero}
\lim_{\alpha\to0} \sup_{ T\in \class}\Pb_{k,\theta}(k <  T < k+N_{\alpha}) =0.
\end{equation}

First, note that 
\[
\alpha \ge \sum_{i=k}^\infty \pi_i \Pb_\infty(T \le i) \ge \Pb_\infty(T \le k) \Pb(\nu \ge k),
\]
and hence,
\begin{equation}\label{Psup}
\hspace{-4mm}\inf_{T\in\class} \Pb_\infty(T > k) \ge 1- \alpha/ \Pb(\nu \ge k), ~~ k \in \Zbb_+,
\end{equation}
which approaches $1$ as $\alpha\to0$ for any fixed $k \in\Zbb_+$. Thus, \eqref{Pinf} follows.

Now, introduce
\[
\begin{aligned} 
U_{\alpha,k}(T) &= e^{(1+\varepsilon) I_\theta N_{\alpha}}\Pb_\infty\brc{k <  T < k+ N_{\alpha}}, 
\\
\beta_{\alpha,k}(\theta)  &= \Pb_{k,\theta}\brc{\frac{1}{N_{\alpha}} \max_{1\le n \le N_{\alpha}} \lambda_{k,k+n}(\theta) \ge (1+\varepsilon) \, I_\theta} .
\end{aligned}
\]
By inequality (3.6) in \cite{TartakovskyVeerTVP05},
\begin{equation}\label{Pktauupper}
 \Pb_{k,\theta}\brc{k <  T < k+ N_{\alpha}} \le  U_{\alpha,k}( T)  + \beta_{\alpha,k}(\theta).
\end{equation}
Using inequality  \eqref{Psup} and the fact that by condition \eqref{Prior}, for all sufficiently large $N_{\alpha}$ (small $\alpha$), there exists a (small) $\delta$ such that
\[
\frac{|\log \Pb(\nu > k+N_{\alpha})|}{ k+N_{\alpha}} \le \mu + \delta,
\]
in just the same way as in the proof of Lemma 1 in \cite{TartakovskyIEEEIT2017} we obtain that
for a sufficiently small $\alpha$
\begin{align} \label{IneqU}
\sup_{ T\in\class} U_{\alpha,k}( T) \le \exp\set{-\frac{I_\theta \varepsilon^2 |\log\alpha|}{I_\theta+\mu + \delta}  + (\mu+\delta) k} .
\end{align}
The right-hand side approaches zero as $\alpha\to0$ for any fixed $k\in \Zbb_+$ and any $\varepsilon>0$ and $\delta>0$.
Also, by condition $\A_1$, $\beta_{\alpha,k}(\theta)\to0$ for all $k\in\Zbb_+$, and therefore, \eqref{Psupclasszero} holds. This completes
 the proof of the lower bound \eqref{LBinclass}.   
\end{IEEEproof}

%%Proof
\begin{IEEEproof}[Proof of Theorem~\ref{Th:FOAOgen}]
(i) For $k \in \Zbb_+$, define the stopping times
\begin{align*}
\tau_A^{(k)}  = & \inf \{n\ge 1: \log \Lambda_{k,k+n}^W +|\log\Pb(\nu\ge k+n)| 
\\
& \quad \ge \log(A/\pi_k)\} .
\end{align*}
Obviously, for any $n >k$,
\begin{align*}
\log S_n^W  & \ge \log\brc{\frac{\pi_k}{\Pb(\nu \ge n)} \Lambda_{k,n}^W} 
\\
&= \log \Lambda_{k,n}^W + \log \pi_k -\log \Pb(\nu \ge n),
\end{align*}
and hence,  for every $A >0$, $(T_A-k)^+ \le  \tau_A^{(k)}$.

Let $N_A=N_{A}(\varepsilon,\theta)=1+\lfloor \log (A/\pi_k) /(I_\theta+\mu-\varepsilon) \rfloor$. Using the same chain of equalities and inequalities as in (A.5) in \cite{TartakovskyIEEEIT2017}, we obtain
that for any $k\in\Zbb_+$ the following inequality holds:
\begin{align}\label{Ektauineq}
& \Eb_{k, \theta}\brcs{(T_A-k)^+}^r   \le \Eb_{k, \theta} \brcs{\brc{\tau_A^{(k)}}^r} \nonumber
\\
 & \le N_{A}^{r} + r 2^{r-1} \sum_{n=N_{A}}^{\infty}  n^{r-1}   \Pb_{k, \theta}\brc{\tau_A^{(k)} >  n}.
 \end{align}
It is easily seen that  for all $k\in\Zbb_+$ and $n\ge N_{A}$
\begin{align*}
&\Pb_{k, \theta}\brc{ \tau_A^{(k)} >n} \nonumber
\\
&\le \Pb_{k, \theta}\Bigg\{\frac{\log \Lambda_{k,k+n}^W}{n}  < \frac{1}{n} \log \brc{\frac{A}{\pi_k}}  \nonumber
\\
& - \frac{|\log\Pb(\nu\ge k+n)|}{n} \Bigg\} \nonumber
\\
& \le \Pb_{k, \theta}\Bigg\{\frac{\log \Lambda_{k, k+n}^W}{n} <I_\theta + \mu -\varepsilon \nonumber
\\
& -\frac{|\log\Pb(\nu\ge k+n)|}{n} \Bigg\} .
\end{align*}
Since, by condition $\mb{CP 1}$, $N_{A}^{-1} |\log\Pb(\nu\ge k+N_{A})| \to \mu$ as $A\to \infty$, for a sufficiently large value of $A$ there exists a small 
$\kappa=\kappa_A$ ($\kappa_A\to 0$ as $A\to\infty$) such that 
\[
\left |\mu - \frac{|\log\Pb(\nu\ge k+N_{A})|}{N_{A}} \right | < \kappa.
\]  
Hence, for all sufficiently large $A$, 
\begin{equation*}
\Pb_{k, \theta}\brc{\tau_A^{(k)} >n} \le \Pb_{k, \theta}\brc{\frac{1}{n}\log \Lambda_{k,k+n}^W <I_\theta -\varepsilon - \kappa} .
\end{equation*}
Also, 
\[
\log  \Lambda_{k, k+n}^W \ge \inf_{\vartheta \in \Gamma_{\delta,\theta}} \lambda_{k, k+n}(\vartheta) +\log W(\Gamma_{\delta,\theta}),
\]
where $\Gamma_{\delta,\theta}=\{\vartheta\in\Theta\,:\,\vert\vartheta-\theta\vert<\delta\}$. Thus, for all sufficiently large  $n$ and $\varepsilon_1 >0$,
\begin{align}
&\Pb_{k, \theta}\brc{ \tau_A^{(k)} >n}   \le \Pb_{k,\theta}\Big(\frac{1}{n}\inf_{\vartheta\in\Gamma_{\delta,\theta}} \lambda_{k, k+n}(\vartheta) \nonumber
\\
&< I_\theta  - \varepsilon- \kappa - \frac{1}{n}\log W(\Gamma_{\delta,\theta})\Big)
\nonumber
\\
&  \le \Pb_{k,\theta}\brc{\frac{1}{n}\inf_{\vartheta\in\Gamma_{\delta,\theta}} \lambda_{k, k+n}(\vartheta) < I_\theta  - \varepsilon_1}. \label{Probktau1}
\end{align}

Using \eqref{Ektauineq}, \eqref{Probktau1}, and inequality $\Pb_\infty(T_A > k) > 1- [A \Pb(\nu >k)]^{-1}$ (see \eqref{Psup}), we obtain
\begin{align}
 & \Rc_{k, \theta}^r(T_A) =\frac{\Eb_{k, \theta}\brcs{\brc{T_A-k}^+}^r}{\Pb_\infty(T_A > k)}  \nonumber
  \\
  &\le   \frac{\brc{1+\left\lfloor \frac{\log (A/\pi_k)}{I_\theta+\mu-\varepsilon} \right\rfloor}^r + r 2^{r-1} \, \Upsilon_{k,r}(\theta, \varepsilon_1)}{1- 1/(A \Pb(\nu \ge k))}. \label{Rckupper}
 \end{align}
Since, by condition $\A_2$, $ \Upsilon_{k,r}(\theta, \varepsilon_1) < \infty$ for all $k\in\Zbb_+$ and $\theta\in\Theta$, this implies the asymptotic upper bound 
\begin{equation}\label{UBkA}
 \Rc_{k, \theta}^m(T_A) \le \brc{\frac{\log A}{I_\theta +\mu}}^m (1+o(1)), \quad A \to \infty
\end{equation}
(for all $0<m \le r$ and $\theta\in\Theta$), which along with the lower bound 
\begin{equation}\label{LBkA}
 \Rc_{k, \theta}^m(T_A) \ge \brc{\frac{\log A}{I_\theta +\mu}}^m (1+o(1)), \quad A \to \infty
\end{equation}
proves the asymptotic relation \eqref{MADDkAOgen}. Note that the lower bound \eqref{LBkA} follows immediately from the lower bound \eqref{LBkinclass} in Lemma~\ref{Lem:LB}
by replacing $\alpha$ with $1/(A+1)$ since it follows from \eqref{PFAMSineq} that $T_A\in \mbb{C}(1/(A+1),\pi)$.

We now get to proving \eqref{MADDAOgen}. Since the MS rule $T_A$ belongs to class $\mbb{C}(1/(A+1), \pi)$,
replacing $\alpha$ by $1/(A+1)$ in the asymptotic lower bound \eqref{LBinclass}, we obtain that under the right-tail condition $\A_1$ the following asymptotic lower bound holds for all $r>0$ and $\theta\in \Theta$:
\begin{equation}\label{LBTA1}
\Rca_{\pi, \theta}^r(T_A) \ge  \brc{\frac{\log A}{I_\theta+\mu}}^r (1+o(1)), ~~ A \to \infty.
\end{equation}
Thus, to prove \eqref{MADDAOgen} it suffices to show that, under the left-tail condition $\A_2$, for $0 <m \le r$ and $\theta\in\Theta$
\begin{equation}\label{UppergenA}
\Rca_{\pi, \theta}^m(T_A) \le  \brc{\frac{\log A}{I_\theta+\mu}}^m (1+o(1), ~~A \to \infty. 
\end{equation}

Using \eqref{Ektauineq} and \eqref{Probktau1}, we obtain that for any $0<\varepsilon < I_\theta+\mu$
\begin{equation}\label{UpperExpr}
\begin{split}
&\Eb^\pi_\theta[(T_{A}-\nu)^+]^r  = \sum_{k=0}^\infty \pi_k  \Eb_{k, \theta}\brcs{(T_A-k)^+}^r
\\
&\le \sum_{k=0}^\infty \pi_k \brc{1+\frac{\log (A/\pi_k)}{I_\theta+\mu-\varepsilon}}^r + 
r 2^{r-1} {\displaystyle\sum_{k=0}^\infty} \pi_k  \Upsilon_{k,r}(\theta,\varepsilon_1).
\end{split}
\end{equation}
This inequality together with the inequality $1-\PFA(T_A) \ge A/(1+A)$ yields
\begin{equation}\label{UpperRca}
\begin{split}
&\Rca_{\pi, \theta}^r(T_A) = \frac{\sum_{k=0}^\infty \pi_k  \Eb_{k,\theta}\brcs{(T_A-k)^+}^r}{1-\PFA(T_A)} 
\\
&\le\frac{{\displaystyle\sum_{k=0}^\infty}  \pi_k \brc{1+\frac{\log (A/\pi_k)}{I_\theta+\mu-\varepsilon}}^r + 
r 2^{r-1} {\displaystyle\sum_{k=0}^\infty} \pi_k  \Upsilon_{k,r}(\theta,\varepsilon_1)}{A/(1+A)} .
\end{split}
\end{equation}
By condition $\A_2$, $\sum_{k=0}^\infty \pi_k  \Upsilon_{k,r}(\theta,\varepsilon_1) < \infty$ for any $\varepsilon_1 >0$ and any $\theta\in\Theta$ 
and, by condition \eqref{Prior1}, $\sum_{k=0}^\infty \pi_k |\log\pi_k|^r < \infty$, which implies that, as $A \to \infty$, for all $0< m \le r$ and all $\theta\in\Theta$
\[
\Rca_{\pi, \theta}^m(T_A) \le \brc{\frac{\log A}{I_\theta+\mu-\varepsilon}}^r (1+o(1)) .
\]
Since $\varepsilon$ can be arbitrarily small, the upper bound \eqref{UppergenA} follows and the proof of the asymptotic expansion  \eqref{MADDAOgen} is complete.

(ii) Setting $A=A_\alpha=(1-\alpha)/\alpha$ in \eqref{MADDkAOgen} and \eqref{MADDAOgen} yields as $\alpha\to 0$ 
\begin{equation}\label{MOapproxgenalpha}
\begin{split}
\Rc^m_{k,\theta}(T_{A_\alpha}) \sim \brc{\frac{|\log\alpha|}{I_\theta+\mu}}^m,
\\
\Rca^m_{\pi,\theta}(T_{A_\alpha}) \sim \brc{\frac{|\log\alpha|}{I_\theta+\mu}}^m,
\end{split}
\end{equation}
which along with the lower bounds \eqref{LBkinclass}  and \eqref{LBinclass}  in Lemma~\ref{Lem:LB} completes the proof of \eqref{FOAOmomentskgen} and \eqref{FOAOmomentsgen}. 
Obviously, all assertions in (ii)  are correct  if threshold $A_\alpha$ is so selected that $T_{A_\alpha}\in\class$ and
$\log A_\alpha\sim |\log \alpha|$ as $\alpha\to0$. The proof is complete.
\end{IEEEproof}

%%Proof
\begin{IEEEproof}[Proof of Lemma~\ref{Lem:LB2}]
Let  $N_{\alpha}= (1-\varepsilon) |\log\alpha|/(I_\theta+\mu_\alpha+\delta_\alpha)$, where $\delta_\alpha>0$ and goes to 0 as $\alpha\to0$. Analogously to \eqref{A1},
\begin{equation} \label{A1new}
\begin{split}
&
\inf_{ T\in\classal}\Rc_{k,\theta}^r(T) 
\\
&\ge N_{\alpha}^r\Big[\inf_{T\in\classal}\Pb_{\infty}( T>k) 
\\
&- \sup_{ T\in \classal}\Pb_{k,\theta}(k <  T < k+N_{\alpha})\Big]
\\
\ge & N_{\alpha}^r\Big[1- \alpha/ \Pb(\nu \ge k) 
\\
&- \sup_{ T\in \classal}\Pb_{k,\theta}(k <  T < k+N_{\alpha})\Big],
\end{split}
\end{equation}
where we used the fact that $\inf_{T\in\classal} \Pb_\infty(T > k) \ge 1- \alpha/ \Pb(\nu \ge k)$ (see \eqref{Psup}). Using \eqref{Pktauupper} and \eqref{IneqU}, we obtain
\begin{align*}
& \sup_{ T\in \classal}\Pb_{k,\theta}(k <  T < k+ N_\alpha)   \le   \beta_{\alpha, k}(\theta)
 \\
 &+ \exp\set{- \frac{\varepsilon^2 I_\theta |\log\alpha| }{I_\theta+\mu_\alpha + \delta_\alpha} + (\mu_\alpha+\delta_\alpha) k}.
 \end{align*}
By condition $\A_1$, $ \beta_{\alpha, k}(\theta)$ goes to zero as $\alpha\to0$ for all $k\in\Zbb_+$. Obviously, the second term vanishes as $\alpha\to0$ for all $\varepsilon\in(0,1)$ and all $k\in\Zbb_+$. 
It follows that, for all $k \in\Zbb_+$ and $\theta\in\Theta$,  
\begin{equation*}
 \sup_{ T\in \classal}\Pb_{k,\theta}(k <  T < k+ N_\alpha)  \to 0 \quad \text{as}~ \alpha\to0
\end{equation*}
and using \eqref{A1new} we obtain that for all $0<\varepsilon<1$,  $r>0$, and $\theta\in\Theta$ as $\alpha\to0$
\[
\inf_{ T\in\classal}\Rc_{\theta,k}^r(T)  \ge (1-\varepsilon)^r\brc{\frac{|\log\alpha|}{I_\theta}}^r (1+o(1)).
\]
Since $\varepsilon$ can be arbitrarily small, the lower bound \eqref{LBinclassk2} follows.
\end{IEEEproof}

%%Proof
\begin{IEEEproof}[Proof of Theorem~\ref{Th:FOAOgen2}]
Setting $A=(1-\alpha)/\alpha$ in inequality \eqref{UpperRca}, we obtain
\begin{align*}
%& 
&\Rca_{\pi^\alpha, \theta}^r(T_A) 
%\\
%& 
\le (1-\alpha)^{-1} \sum_{k=0}^\infty  \pi_k^\alpha \brc{1+\frac{\log ((1-\alpha)/\alpha\pi_k^\alpha)}{I_\theta+\mu_\alpha-\varepsilon}}^r 
\\
& + 
r 2^{r-1} \sup_{k\in \Zbb_+}\Upsilon_{k,r}(\theta,\varepsilon_1) .
\end{align*}
Using conditions \eqref{Prior3} and $\A_3$ and taking into account that $\mu_\alpha\to0$ as $\alpha\to0$ yields
\[
\Rca_{\pi^\alpha, \theta}^r(T_{A_\alpha})   \le  \brc{\frac{|\log\alpha|}{I_\theta-\varepsilon}}^r (1+o(1)).
\]
Since $\varepsilon$ can be arbitrary small,  we obtain the asymptotic upper bound 
\[
\Rca_{\pi^\alpha, \theta}^r(T_{A_\alpha})  \le  \brc{\frac{|\log\alpha|}{I_\theta}}^r (1+o(1)),
\]
as $\alpha \to 0$, which along with the lower bound \eqref{LBinclass2}  proves \eqref{MADDAOgen2}. Clearly, this upper bound holds if we take any $A=A_\alpha$ such that $\PFA(T_{A_\alpha}) \le \alpha$ and 
$\log A_\alpha \sim |\log \alpha|$ as $\alpha \to 0$.

Next, substituting  $A=(1-\alpha)/\alpha$  in \eqref{Rckupper} (or more generally any $A$ such that $\log A \sim |\log \alpha|$ and $\PFA(T_A) \le \alpha$), we obtain
\begin{align*}
 & \Rc_{k, \theta}^r(T_A)
\le   
\\
&\frac{\brcs{1+\left\lfloor \log \frac{((\alpha/(1-\alpha)\pi_k^\alpha)}{I_\theta+\mu_\alpha-\varepsilon} \right \rfloor}^r + r 2^{r-1} \, \Upsilon_{k,r}(\theta, \varepsilon_1+ \delta_\alpha)}{1- 1/(A \Pb(\nu \ge k))},
 \end{align*}
which due to conditions $\A_3$ and \eqref{Prior3} and the fact that $\mu_\alpha, \delta_\alpha \to 0$ implies that, for all fixed $k \in \Zbb_+$ and all $\theta\in\Theta$ as $\alpha\to0$,
\[
\Rc_{k, \theta}^r(T_A) \le \brc{\frac{|\log\alpha|}{I_\theta}}^r (1+o(1)).
\]
This upper bound together with the lower bound \eqref{LBinclassk2} yields  \eqref{MADDAOgenk2} and the proof is complete.
\end{IEEEproof}

%%Proof
\begin{IEEEproof}[Proof of Theorem~\ref{Th:AOCMSR}]
(i)   For $\varepsilon\in(0,1)$, let $M_A=M_{A}(\varepsilon,\theta) = (1-\varepsilon) I_\theta^{-1} \log A$. Similarly to \eqref{A1} we obtain
\begin{equation} \label{LBSRA}
\begin{split}
\Rc_{k,\theta}^r(\wtT_A)  &\ge M_{A}^r\big[\Pb_{\infty}(\wtT_A>k) 
\\
&- \Pb_{k,\theta}(k <  \wtT_A < k+M_{A})\big]
\end{split}
\end{equation}
and similarly to \eqref{Pktauupper},
\begin{equation}\label{PkTBupper}
 \Pb_{k,\theta}\brc{0 < \wtT_A -k < M_{A}} \le  U_{A,k}(\wtT_A)  + \beta_{A,k}(\theta) ,
\end{equation}
where
\[
\begin{aligned} 
U_{A,k}(\wtT_A) & = e^{(1+\varepsilon) I_\theta M_{A}}\Pb_\infty\brc{0 < \wtT_A - k <M_{A}}, 
\\
\beta_{A,k}(\theta)  &= \Pb_{k,\theta}\brc{\frac{1}{M_{A}} \max_{1\le n \le M_{A}} \lambda_{k,k+n}(\theta) \ge (1+\varepsilon) \, I_\theta} .
\end{aligned}
\]
 Since 
\[
\begin{aligned}
\Pb_\infty\brc{0 < \wtT_A - k <M_{A}} & \le \Pb_\infty\brc{\wtT_A  < k+ M_{A}} 
\\
&\le (k+\omega + M_{A})/A,
\end{aligned}
\]
we have
\begin{equation}\label{UpperU}
 U_{A,k}(\wtT_A) \le \frac{k+ \omega+(1-\varepsilon) I_\theta^{-1} \log A}{A^{\varepsilon^2}}.
\end{equation}
Therefore,  $U_{A,k}(\wtT_A)\to 0$ as $A\to\infty$ for any fixed $k$. Also, $\beta_{A,k}(\theta)\to 0$ by condition $\C_1$, so that $\Pb_k\brc{0 < \wtT_A -k < M_{A}}\to0$ for any fixed $k$. Since
$\Pb_{\infty}(\wtT_A>k) > 1- (\omega+k)/A$, it follows from \eqref{LBSRA}  that for an arbitrary $\varepsilon \in (0,1)$ as $A\to\infty$
\begin{equation*}
\Rc_{k,\theta}^r(\wtT_A) \ge \brc{\frac{(1-\varepsilon) \log A}{I_\theta}}^r (1+o(1)),
\end{equation*}
which yields the asymptotic lower bound (for any fixed $k\in\Zbb_+$ and $\theta\in\Theta$)
\begin{equation}\label{LBSRAasympt}
\Rc_{k,\theta}^r(\wtT_A) \ge \brc{\frac{ \log A}{I_\theta}}^r (1+o(1)), \quad A\to\infty.
\end{equation}

To prove \eqref{MomentskMSR1} it suffices to show that
\begin{equation}\label{UBSRAasympt}
\Rc_{k,\theta}^r(\wtT_A) \le \brc{\frac{ \log A}{I_\theta}}^r (1+o(1)), \quad A\to\infty.
\end{equation}

For $k \in \Zbb_+$, define the stopping times
\begin{align*}
\tilde{\tau}_A^{(k)}  = \inf \{n\ge 1: \log \Lambda_{k,k+n}^W \ge \log A\} .
\end{align*}
Obviously, for any $n >k$, $\log R_n^W  \ge \log \Lambda_{k,n}^W$, and hence,  for every $A >0$, $(\wtT_A-k)^+ \le  \tilde{\tau}_A^{(k)}$.
Analogously to \eqref{Ektauineq}, we have
\begin{align}\label{Ektildetauineq}
& \Eb_{k, \theta}\brcs{(\wtT_A-k)^+}^r   \le \Eb_{k, \theta} \brcs{\brc{\tilde{\tau}_A^{(k)}}^r} \nonumber
\\
 & \le \widetilde{M}_{A}^{r} + r 2^{r-1} \sum_{n= \widetilde{M}_{A}}^{\infty}  n^{r-1}   \Pb_{k, \theta}\brc{\tilde{\tau}_A^{(k)} >  n},
 \end{align}
where $\widetilde{M}_{A}= \widetilde{M}_{A}(\varepsilon,\theta)=1+\lfloor \log (A)/(I_\theta-\varepsilon) \rfloor$. For a sufficiently large $n$ similarly to \eqref{Probktau1} we have
\begin{align}
\Pb_{k, \theta}\brc{ \tilde\tau_A^{(k)} >n}  \le \Pb_{k,\theta}\brc{\frac{1}{n}\inf_{\vartheta\in\Gamma_{\delta,\theta}} \lambda_{k, k+n}(\vartheta) < I_\theta  -\varepsilon}. \label{Probktildetau1}
\end{align}
Using \eqref{Ektildetauineq} and \eqref{Probktildetau1}, we obtain the inequality
\begin{equation}\label{Ekinequpper}
\begin{split}
&\Eb_{k, \theta}\brcs{\brc{\wtT_A-k}^+}^r 
\\
& \le \brc{1+\left\lfloor \frac{\log A} {I_\theta-\varepsilon} \right\rfloor}^r + r 2^{r-1} \, \Upsilon_{k,r}(\theta, \varepsilon),
\end{split}
\end{equation}
which along with the inequality $\Pb_\infty(\wtT_A > k) > 1- (\omega +k)/A$ implies the inequality
\begin{align}
 & \Rc_{k, \theta}^r(\wtT_A) =\frac{\Eb_{k, \theta}\brcs{\brc{\wtT_A-k}^+}^r}{\Pb_\infty(\wtT_A > k)}  \nonumber
  \\
  &\le   \frac{\brc{1+\left\lfloor \frac{\log A}{I_\theta-\varepsilon} \right\rfloor}^r + r 2^{r-1} \, \Upsilon_{k,r}(\theta, \varepsilon)}{1- (\omega +k)/A}. \label{RckupperSR}
 \end{align}
Since, by condition $\A_2$, $ \Upsilon_{k,r}(\theta, \varepsilon) < \infty$ for all $k\in\Zbb_+$ and $\theta\in\Theta$, this implies the asymptotic upper bound \eqref{UBSRAasympt} 
and completes the proof of the asymptotic approximation  \eqref{MomentskMSR1}.

We now continue with proving  \eqref{MomentsMSR1}. By the Chebyshev inequality,
\begin{align}
& \Rca_{\pi, \theta}^r(\wtT_A) \ge \Eb_{\theta}^\pi[( \wtT_A-k)^+]^r\nonumber
\\
&\ge M_{A}^r \Pb_{\theta}^\pi(\wtT_A  -\nu > M_{A})\nonumber
\\
& \ge  M_{A}^r\brcs{\Pb_{\theta}^\pi(\wtT_A> \nu)-\Pb_{\theta}^\pi(\nu <  T < \nu+ M_{A})}  \nonumber
\\
& \ge M_{A}^r\brcs{1-\frac{\bar\nu+\omega}{A} - \Pb^\pi_\theta \brc{0 < \wtT_A-\nu < M_{A}}} . \label{LBMA}
\end{align}
Let $K_A$ be an integer number that approaches infinity as $A\to\infty$. Using \eqref{PkTBupper} and \eqref{UpperU}, we obtain the following upper bound 
\begin{align}\label{ProbineqSR}
&\Pb^\pi_\theta(0< \wtT_A -\nu < M_{A})  \nonumber
\\
&= \sum_{k=0}^\infty  \pi_k  \Pb_{k,\theta}\brc{0 < \wtT_A -k < M_{A}} \nonumber
 \\
 &\le \Pb(\nu > K_{A}) +   \sum_{k=0}^\infty \pi_k U_{A,k}(\wtT_A) + \sum_{k=0}^{K_{A}}  \pi_k  \beta_{A, k} \nonumber
 \\ 
 & \le  \Pb(\nu > K_{A}) \nonumber
 \\
 & + \frac{\bar\nu+\omega+ (1-\varepsilon) I_\theta^{-1} \log A}{A^{\varepsilon^2}} 
  +  \sum_{k=0}^{K_{A}}  \pi_k  \beta_{A,k},
 \end{align}
where the first two  terms go to zero as $A\to\infty$ since $\bar\nu$ and $\omega$ are finite (by Markov's inequality $\Pb(\nu>K_A) \le \bar\nu/K_A$) 
and the last term also goes to zero  
by condition $\A_1$ and Lebesgue's dominated convergence theorem. Thus, for all $0<\varepsilon <1$,
$\Pb^\pi_\theta(0< \wtT_A -\nu < M_{A})$ approaches $0$ as $A\to\infty$.
Using inequality \eqref{LBMA}, we obtain that for any $0<\varepsilon <1$ as $A\to\infty$
\[
\Rca_{\pi,\theta}^r(\wtT_A) \ge (1-\varepsilon)^r \brc{\frac{\log A}{I_\theta}}^r (1+o(1)),
\]
which yields the asymptotic lower bound (for any $r>0$ and $\theta\in\Theta$)
\begin{equation}\label{LBA}
\hspace{-2mm}\Rca_{\pi,\theta}^r(\wtT_A) \ge \brc{\frac{\log A}{I_\theta}}^r (1+o(1)), \quad A \to \infty.
\end{equation}

To obtain the upper bound it suffices to use inequality \eqref{Ekinequpper}, which along with the fact that $\Pb^\pi(\wtT_A > \nu) > 1- (\bar\nu +\omega)/A$ yields (for every $0<\varepsilon < I_\theta$)
\begin{align}\label{RcarMSR}
&\Rca_{\pi,\theta}^r(\wtT_A) = \frac{\sum_{k=0}^\infty \pi_k  \Eb_{k}[(\wtT_A-k)^+]^r}{\Pb^\pi(\wtT_A > \nu)} \nonumber
\\
&\le\frac{\brc{1+\frac{\log A}{I_\theta-\varepsilon}}^r + r 2^{r-1}\, \sum_{k=0}^\infty \pi_k  \Upsilon_{k,r}(\theta,\varepsilon)}{1-(\omega+\bar\nu)/A} .
\end{align}
Since by condition $\A_2$,
\[
 \sum_{k=0}^\infty \pi_k  \Upsilon_{k,r}(\theta,\varepsilon) < \infty \quad \text{for any}~ \varepsilon >0~\text{and}~\theta\in\Theta,
\]
we obtain that, for every $0<\varepsilon < I_\theta$ as $A \to \infty$,
\[
\Rca_{\pi,\theta}^r(\wtT_A) \le \brc{\frac{\log A}{I_\theta-\varepsilon}}^r (1+o(1)) .
\]
Since $\varepsilon$ can be arbitrarily small, this implies the asymptotic (as $A\to\infty$) upper bound
\begin{equation}\label{UBB}
\Rca_{\pi,\theta}^r(\wtT_A) \le \brc{\frac{\log A}{I_\theta}}^r (1+o(1)),
\end{equation}
which along with the lower bound \eqref{LBA} completes the proof of (i).

(ii) To prove \eqref{MomentsMSR2} and \eqref{MomentskMSR2} it suffices to substitute $\log A \sim |\log \alpha|$ (in particular, we may take $A=(b \, \omega +\bar{\nu})/\alpha$) in 
\eqref{MomentsMSR1} and \eqref{MomentskMSR1}.
\end{IEEEproof}

%%Proof
\begin{IEEEproof}[Proof of Theorem~\ref{Th:AoptMSR}]
Previous results make the proof elementary. Indeed, substitution $A=A_\alpha=(b_\alpha\omega_\alpha + \bar\nu_\alpha)/\alpha$ in \eqref{RcarMSR} yields the upper bound
\begin{align*}
& \Rca_{\pi^\alpha,\theta}^r(\wtT_{A_\alpha}) \le
\\
& \frac{\brc{1+\frac{\log ((\omega_\alpha + \bar\nu_\alpha)/\alpha)}{I_\theta-\varepsilon}}^r + r 2^{r-1}\, \sup_{k\ge 0} \Upsilon_{k,r}(\theta,\varepsilon)}{1-\alpha} ,
\end{align*}
which implies the asymptotic upper bound
\[
\Rca_{\pi^\alpha,\theta}^r(\wtT_{A_\alpha}) \le \brc{\frac{|\log\alpha|}{I_\theta}}^r (1+o(1)), \quad \alpha \to 0,
\]
since, by condition \eqref{Prior4},  $\log [(\omega_\alpha + \bar\nu_\alpha)/\alpha] \sim |\log\alpha|$ and by condition $\A_3$, $\sup_{k\ge 0} \Upsilon_{k,r}(\theta,\varepsilon)<\infty$ for all $\theta\in\Theta$.
This upper bound along with the lower bound \eqref{LBinclass2} in Lemma~\ref{Lem:LB2} proves \eqref{MomentsMSRmu0}.

Finally, the asymptotic upper bound
\begin{align*}
\Rc_{k,\theta}^r(\wtT_{A_\alpha}) 
\le \brc{\frac{|\log\alpha|}{I_\theta}}^r (1+o(1)), \quad \alpha \to 0,
\end{align*}
for all $k \in \Zbb_+$ and $\theta\in\Theta$ follows immediately from \eqref{RckupperSR} with $A=A_\alpha=(c_\alpha\omega_\alpha + \bar\nu_\alpha)/\alpha$, 
which along with the lower bound \eqref{LBinclassk2} in Lemma~\ref{Lem:LB2} proves \eqref{MomentskMSRmu0}.
\end{IEEEproof}
%%%

%%Proof
\begin{IEEEproof}[Proof of Theorem~\ref{Th:FOasopt_pureBayes}]
Since $\Theta$ is compact it follows from the asymptotic approximation \eqref{MADDAOgen} in Theorem~\ref{Th:FOAOgen} that as $A\to\infty$
\[
\begin{split}
\int_\Theta \Rca_{\pi, \vartheta}^r(T_{A}) \, \drm W(\vartheta)  &=  \int_\Theta \brc{\frac{\log  A}{I_\vartheta+\mu}}^r  \drm W(\vartheta) (1+o(1))
\\
& = D_{\mu,r} (\log A)^r (1+o(1)) .
\end{split}
\]
Since $\PFA(T_A) < 1/A$, we obtain the following asymptotic approximation for the integrated risk 
\begin{equation}\label{FOAapprIR}
\rho_{\pi, W}^{c,r}(T_A) \sim D_{\mu,r} \, c \, (\log A)^r \quad \text{as}~ A \to \infty.
\end{equation}

Next, it is easily seen that, for any $r \ge 1$ and $\mu \in (0,\infty)$, threshold~$A_{c,r}$ goes to
infinity as $c\to 0$ with such rate  that $\log A_{c,r} \sim |\log c|$. As a result, we obtain that  
\[
\rho_{\pi, W}^{c,r}(T_{A_{c,r}}) \sim D_{\mu,r} \, c \, |\log c|^r  \quad \text{as}~~c\to 0.
\]

All it remains to do is to prove the lower bound
\begin{equation} \label{LBIR}
\inf_{T\ge 0}\rho_{\pi, W}^{c,r}(T)\ge D_{\mu,r} \, c \, |\log c|^r (1+o(1)) \quad \text{as}~~ c \to 0.
\end{equation}
In fact, since
\[
\lim_{c\to 0} \frac{G_{c,r}(A_{c,r})}{D_{\mu,r} \, c \, |\log c|^r} =1,
\]
it suffices to prove that
\begin{equation}\label{needtoprove}
\frac{\inf_{T\ge 0} \rho_{\pi, W}^{c,r}(T)}{G_{c,r}(A_{c,r})} \ge 1 +o(1) \quad \text{as}~ c \to 0.
\end{equation}
This can be done by contradiction. Indeed, suppose that \eqref{needtoprove} is wrong, i.e., there exists a stopping rule $T=T_{c}$ such that
\begin{equation} \label{assumptionTc}
\frac{\rho_{\pi, W}^{c,r}(T_c)}{G_{c,r}(A_{c,r})} < 1 +o(1) \quad \text{as}~ c \to 0 .
\end{equation}
Let $\alpha_c=\PFA(T_c)$. First, $\alpha_c \to 0$ as $c\to 0$ since
\[
\alpha_c\le \rho_{c,r}^{\pi,W}(T_c) < G_{c,r}(A_{c,r})(1+o(1))\to 0 \quad \text{as}~~c\to0.
\]
Second, it follows from Lemma~\ref{Lem:LB} that, as $\alpha_c \to 0$, 
\[
\Rc_{\pi,W}^r(T_c)  \ge \int_\Theta (I_\vartheta + \mu)^{-r} \drm W(\vartheta) |\log \alpha_c|^r (1+o(1)),
\]
and hence, as $c \to0$,
\begin{align*}
\rho_{\pi, W}^{c,r}(T_c)& = \alpha_c + c\, (1-\alpha_c) \Rc_{\pi,W}^r(T_c) 
\\
& \ge \alpha_c + c \, D_{\mu,r}  |\log \alpha_c|^r (1+o(1)) .
\end{align*}
Thus,
\begin{align*}
\frac{\rho_{\pi, W}^{c,r}(T_c)}{G_{c,r}(A_{c,r})} &\ge \frac{G_{c,r}(1/\alpha_c) +c \, D_{\mu,r} |\log \alpha_c|^r o(1)}{\min_{A>0} G_{c,r}(A)} 
\\
&\ge 1+o(1),
\end{align*}
which contradicts~\eqref{assumptionTc}. Hence, \eqref{needtoprove} follows and the proof is complete.
\end{IEEEproof}

%%Proof
\begin{IEEEproof}[Proof of Theorem~\ref{Th:FOaopt_pureBayesflat}]  
(i) Using \eqref{UpperExpr}, we obtain
\[
\begin{split}
&\int_\Theta \Eb^{\pi^c}_\vartheta[(T_{A}-\nu)^+]^r \, \drm W(\vartheta) 
\\
& \le \sum_{k=0}^\infty \pi_k^c\int_\Theta \brc{1+\frac{\log (A/\pi_k^c)}{I_\vartheta+\mu_c-\varepsilon}}^r  \drm W(\vartheta) 
\\
& + r 2^{r-1} \int_\Theta \Upsilon_{r}(\vartheta,\varepsilon_1) \, \drm W(\vartheta),
\end{split}
\]
where the last term is finite since  $\Theta$ is compact and $\Upsilon_r(\vartheta,\varepsilon_1) < \infty$ for all $\vartheta\in\Theta$
due to condition $\A_3$ and where 
\[
\sum_{k=0}^\infty  \pi_k^c |\log \pi_k^c|^r =o(|\log c|^r) \quad \text{as}~ c \to 0
\]
by assumption \eqref{Prior3c}. Recall that, as established in the proof of Theorem~\ref{Th:FOasopt_pureBayes} above, $\log A_{c,r} \sim |\log c|$ as $c\to0$ when $A_{c,r}$ satisfies \eqref{threshold1} 
and that $\mu_c\to0$. Therefore,  as $c\to0$,
\begin{align*}
&\int_\Theta \Eb^{\pi^c}_\vartheta[(T_{A_{c,r}}-\nu)^+]^r \, \drm W(\vartheta) 
\\
&\le  \int_\Theta \brc{\frac{|\log c|}{I_\vartheta-\varepsilon}}^r \drm W(\vartheta) (1+o(1)) .
\end{align*}
Since $\varepsilon$ can be arbitrary small, we obtain 
\begin{align*}
&c \, \int_\Theta \Eb^{\pi^c}_\vartheta[(T_{A_{c,r}}-\nu)^+]^r \, \drm W(\vartheta) 
\\
&\le D_{r} \, c \,  |\log c|^r(1+o(1)) \quad \text{as}~c\to0.
\end{align*}
Since $\PFA(T_{A_{c,r}}) < 1/A_{c,r}= o(c|\log c|^r)$ it follows that   
\begin{equation}\label{UpperMSIRflat}
\rho_{\pi^c, W}^{c,r}(T_{A_{c,r}}) \le D_{r} \, c \, |\log c|^r (1+o(1)) \quad \text{as}~c\to0.
\end{equation}

The lower bound
\begin{equation} \label{LBIR2}
\inf_{T\ge 0}\rho_{\pi^c, W}^{c,r}(T)\ge D_{r} \, c \, |\log c|^r (1+o(1)) \quad \text{as}~~ c \to 0
\end{equation}
can be deduced using Lemma~\ref{Lem:LB2} and the argument essentially similar to that used in the proof of the lower bound \eqref{LBIR} above with $G_{c,r}(A)=1/A + c \, D_{r} (\log A)^r$.

Using the asymptotic upper bound \eqref{UpperMSIRflat} and the lower bound \eqref{LBIR2} simultaneously, we obtain \eqref{MomentsMSmu0Bayesflat}, which completes the proof of (i).

(ii) In order to prove \eqref{MomentsMSRmu0Bayesflat} it suffices to prove the asymptotic upper bound  
\begin{equation}\label{UpperMSRIRflat}
\rho_{\pi^c, W}^{c,r}(\wtT_{A_{c,r}}) \le D_{r} \, c \, |\log c|^r (1+o(1)) \quad \text{as}~c\to0.
\end{equation}
Define 
\[
\widetilde{G}_{c,r}(A)=(\omega_c b_c +\bar{\nu}_c)/A + c \, D_{r} (\log A)^r.
\]
Threshold $A_{c,r}$ that satisfies equation \eqref{threshold2} minimizes $\widetilde{G}_{c,r}(A)$, and it is easily seen that $\log A_{c,r} \sim |\log c|$ as $c\to 0$ 
since,  by assumption \eqref{Prior4c},  $\omega_c b_c +\bar{\nu}_c = o(|\log c|)$ as $c\to0$.

Using inequality \eqref{Ekinequpper}, we obtain
\begin{align*}
& \int_\Theta \Eb^{\pi^c}_\vartheta[(\wtT_{A}-\nu)^+]^r \, \drm W(\vartheta) 
\\
& \le \int_\Theta \brc{1+\frac{\log A}{I_\vartheta-\varepsilon}}^r \, \drm W(\vartheta)
\\ 
& \quad + r 2^{r-1} \int_\Theta \Upsilon_{r}(\vartheta,\varepsilon) \, \drm W(\vartheta) ,
\end{align*}
where the last term is finite since  $\Theta$ is compact and $\Upsilon_r(\vartheta,\varepsilon) < \infty$ for all $\vartheta\in\Theta$
due to condition $\A_3$.  Therefore, for an arbitrary small $\varepsilon$ as $c \to 0$,
\begin{align*}
&\int_\Theta \Eb^{\pi^c}_\vartheta[(\wtT_{A}-\nu)^+]^r \, \drm W(\vartheta)  
\\
&\le \int_\Theta \brc{\frac{|\log c|}{I_\vartheta-\varepsilon}}^r \, \drm W(\vartheta) (1+o(1)),
\end{align*}
which implies that, as $c\to0$, 
\[
c \, \int_\Theta \Eb^{\pi^c}_\vartheta[(\wtT_{A}-\nu)^+]^r \, \drm W(\vartheta)  \le D_r \, c \, |\log c|^r (1+o(1)).
\]
Since $\PFA(\wtT_{A_{c,r}}) \le (\omega_c b_c+\bar\nu_c)/A_{c,r} = o(c|\log c|^r)$, we obtain \eqref{UpperMSRIRflat}, which along with the lower bound \eqref{LBIR2} proves \eqref{MomentsMSRmu0Bayesflat}.
The proof of (ii) is complete.
\end{IEEEproof}

%%%%%%%%%%%%%%%%%%%%%%%%%%%%%%%%%%%%%%%%%%%%%%
%\bibliographystyle{IEEEtran}

%\bibliography{\string~/Library/Mobile Documents/SashaDocs/Papers(my)/Bibliography/Bib_main}
%\bibliography{~/Library/Mobile Documents/Bibliography/Bib_main.bib}
%\bibliography{../../Bib_main.bib}
%\bibliography{Bib_main}


%%%%% Biography%%%%%%%%
\begin{IEEEbiography}[{\includegraphics[width=1.2in,height=1.45in,keepaspectratio]{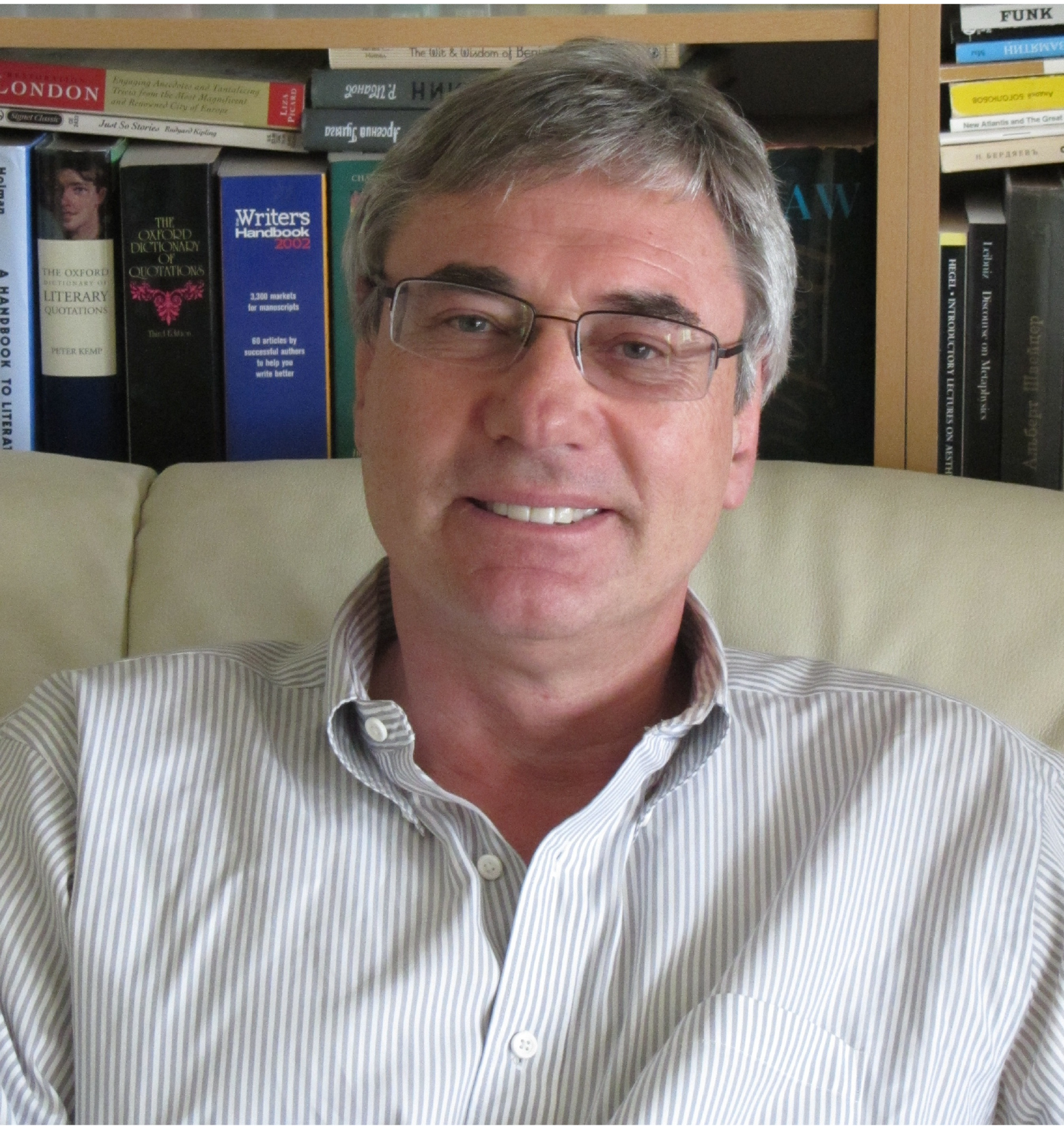}}]{Alexander G.\ Tartakovsky} (M'01-SM'02)  
research interests include theoretical and applied statistics; applied
probability; sequential analysis; changepoint detection phenomena; and a variety of
applications including statistical image and signal processing; video tracking; detection
and tracking of targets in radar and infrared search and track systems; near-Earth space informatics; information
integration/fusion; intrusion detection and network security; and detection and tracking
of malicious activity. He is the author of two books (the third is in preparation),
several book chapters, and over 100 papers. Dr.\ Tartakovsky is a {\em Fellow of the Institute of Mathematical Statistics} and a senior
member of IEEE. He received several awards, including a {\em 2007 Abraham Wald Award in
Sequential Analysis.} Dr.\ Tartakovsky obtained a Ph.D.\ degree and an
advanced Doctor-of-Science degree both from Moscow Institute of
Physics and Technology, Russia (FizTech). During 1981--92, he was first a Senior Research Scientist and then a
Department Head at the Institute of Radio Technology (Moscow, Russian Academy of
Sciences) as well as a Professor at FizTech, working on the application of statistical methods to optimization and modeling
of information systems. From 1993 to 1996, Dr.\ Tartakovsky worked at the University
of California, Los Angeles (UCLA), first in the Department of Electrical Engineering and then
in the Department of Mathematics. From 1997 to 2013, he was a Professor in the Department of Mathematics
and the Associate Director of the Center
for Applied Mathematical Sciences at the University of Southern California (USC), Los
Angeles. From 2013 to 2015, he was a Professor of Statistics in the Department of Statistics at the University of Connecticut, Storrs. 
Currently, he is the Head of the Space Informatics Laboratory at FizTech as well as Vice President of  AGT StatConsult, Los Angeles, California.
\end{IEEEbiography}

\end{document}